\documentclass[12pt]{article}

\usepackage[T2A]{fontenc}
\usepackage[margin=1in]{geometry} 
\usepackage{amsmath,amsthm,amssymb}
\usepackage[utf8]{inputenc}
\usepackage[english]{babel}
\usepackage{graphicx}
\usepackage{tablefootnote}
\usepackage{float}

\usepackage{amsfonts, amsmath, amsthm}
\usepackage{amssymb}
\usepackage{url}
\usepackage{cite}

\usepackage{indentfirst}

\usepackage{mathtools}
\usepackage{makecell}
\usepackage{tikz}
\usepackage{dsfont}
\usepackage{secdot}
\usepackage{graphicx}

\usepackage{algorithm}
\usepackage{algorithmic}
\usepackage{float}
\usepackage{comment}

\usepackage{multicol}
\usepackage{lipsum}
\usepackage{mwe}

\usepackage{caption}
\usepackage{subcaption}

\usepackage[unicode=true]{hyperref}

\sectiondot{section}
\sectiondot{subsection}
\theoremstyle{definition}

\newtheorem{theorem}{Theorem}
\newtheorem{lemma}{Lemma}

\newtheorem{definition}{Definition}
\newtheorem{remark}{Remark}
\newtheorem{statement}{Proposition}  

\newtheorem{question}{Question}
\makeatletter
\renewcommand{\ALG@name}{Algorithm}
\newcommand{\sgn}{\operatorname{sgn}}

\allowdisplaybreaks

\graphicspath{{images/}}

\begin{document}

\begin{center}
\textbf{Optimal partitions of the flat torus into parts of smaller diameter}

{\bf
D.\,S.\,Protasov$^{1,3*}$,
A.\,D.\,Tolmachev$^{1,2,3**}$,
V.\,A.\,Voronov$^{1,3,***}$
}

\textit{$^1$Moscow Institute of Physics and Technology, Dolgoprudniy, Russia} \\
\textit{$^2$Skolkovo Institute of Science and Technology, Moscow, Russia} \\
\textit{$^3$Caucasus Mathematical Center of Adyghe State University, Maikop, Russia}

$^{*}$ E-mail: dmitry.protasov@gmail.com\\
$^{**}$ E-mail: tolmachev.ad@phystech.edu\\
$^{***}$ E-mail: v-vor@yandex.ru\\

\renewcommand{\abstractname}{\vspace{-\baselineskip}}
\begin{abstract}
We consider the problem of partitioning a two-dimensional flat torus $T^2$ into $m$ sets in order to minimize the maximal diameter of a part. For $m \leqslant 25$ we give numerical estimates for the maximal diameter $d_m(T^2)$ at which the partition exists. Several approaches are proposed to obtain such estimates. In particular, we use the search for mesh partitions via the SAT solver, the global optimization approach for polygonal partitions, and the optimization of periodic hexagonal tilings. For $m=3$, the exact estimate is proved using elementary topological reasoning.
\end{abstract}
\end{center}

\section{Introduction}\label{introduction}

The present paper is devoted to one of the special questions associated with the Borsuk problem. Recall that the main statement of the Borsuk problem was that it is possible to divide any bounded subset of $n$-dimensional Euclidean space into $n+1$ parts of strictly smaller diameter \cite{borsuk,borsuk1}. This statement is known to be false for $n\geqslant 64$, and for $4 \leqslant n \leqslant 63$ it is not known whether it is true  \cite{bondarenko2014borsuk,jenrich201464}. 

Moving on to a particular problem, we consider all possible partitions of a given bounded set into $m$ parts. Suppose that the diameter of the parts is bounded by some number $d$. What is a minimal $d$ for which the partition exists? Unlike numerous packing and covering problems, this question is rather poorly studied. Apparently, the first statements of this kind belong to Heppes\cite{heppes}. More attention has been paid to the maximal diameter of parts, at which it will be possible to divide any set $F \subset \mathbb{R}^n$ of unit diameter \cite{Lassak,Filimonov,borsuk2}. Recent work  has devoted to a similar problem for $\mathbb{R}^n$ with an arbitrary norm \cite{lian2021partition}.
 
Let $(X, \rho)$ be a metric space, let $F \subset X$ be a bounded set. The diameter of the set $F$ is defined as follows: $$ \operatorname{diam} F  = \sup_{x, y \in F} \rho(x, y) $$ 

Next, for $m \geqslant 1$, we define the following quantity.
$$d_m(F) = \inf\{x \in \mathbb{R}^{+}\; |\; \exists F_{1}, \ldots, F_{m} \subset X: F = F_{1} \sqcup \ldots \sqcup F_{m}, \; \forall i \:  \operatorname{diam}(F_{i}) \leqslant x \}.$$

Note that the problem of partitioning the set $F$ into arbitrary sets $F_i$ in order to minimize the maximal diameter of $F_i$, $i = 1,2, \dots, m$ is obviously equivalent to the problem of covering by arbitrary sets. In addition, one can assume that all sets are closed (or all are open).
$$d_m(F) = \inf\{x \in \mathbb{R}^{+} \;|\; \exists F_{1}, \ldots, F_{m} \subset X: F \subseteq F_{1} \cup \ldots \cup F_{m}, \; \forall i \:  \operatorname{diam}(F_{i}) \leqslant x \}.$$
 
 
 For an Euclidean metric in $\mathbb{R}^n$ and $m = 1,2, \dots$ one can define the value $d_{n,m} = \sup d_{m}(F)$  where the suprema are taken over all sets \(F \subset\mathbb{R}^n\) of unit diameter \cite{borsuk2, borsuk3, bikeev2023}. 
 
 In \cite{Lenz,Lassak,Filimonov,koval2023partition} several estimates of $d_{2,m}, m=3, 4 \dots$ have been found. In  \cite{DAM_our_article} some upper and lower bounds for $d_{2,m}$ and $d_{3,4}$ was improved. 

In this paper, we consider the problem for a particular set, namely a two-dimensional flat torus. Compared to polygon partitioning \cite{heppes,voronov2023_motor}, this problem is in some sense more interesting. 
 In the case of a flat torus the problem contains no geometric constraints. As a consequence, it is impossible to prove an exact lower estimate by considering a finite set of points, since the optimal partition can be shifted in such a way that these points do not belong to the boundaries. 

Consider the flat torus $T^2$ as the factor-space $T^2 =\mathbb{R}^2/\mathbb{Z}^2$. The corresponding metric $\rho_T$ is obtained from the Euclidean metric on the plane. 
Namely, we define the metric $\rho_T$ on the flat torus as follows: 
\[\rho_T\left(u,v\right)  = \left[\left(\min\{|x_1 - x_2|, 1 - |x_1 - x_2|\}\right)^2 + \left(\min\{|y_1 - y_2|, 1 - |y_1 - y_2|\}\right)^2\right]^{1/2},
\] 
where $u = (x_1,y_1), v = (x_2, y_2)$.

In the following, we will consider the estimates of $d_m(T^2)$.
Note that in any dimension it makes no sense to ask a question about the minimal number of parts at which the diameter of the part is less than the diameter of the torus. Indeed, a flat torus of arbitrary dimension $T^n$ can obviously be divided into three layers of thickness $\frac{1}{3}$, which will have diameter $$\left(\frac{n-1}{4}+\frac{1}{9}\right)^{1/2}<\frac{n^{1/2}}{2}=\operatorname{diam} T^n.$$ In other words, the Borsuk number of the flat torus $T^n$ is equal to $3$ for all $n \geqslant 1$.  

\section{Main result}

First, we give estimates based on simple geometric considerations.
 
 \theorem \label{theorem1} The following inequalities are satisfied:
  
 \begin{equation} d_m(T^2) \leqslant \sqrt{\frac{1}{4} + \frac{1}{m^2}} , \quad m = 1,2, \dots  \label{tor_stripes} \end{equation}
 
 \begin {equation} d_m(T^2) \geqslant \frac{2}{\sqrt{\pi m}} , \quad m \geqslant 6 \label{tor_sqrt} \end {equation}

\begin {equation} d_{k^2 + k - 1}(T^2) \geqslant \frac{1}{k}, \quad k =1, 2, \dots  \label{tor_lower_lines} \end {equation}

Note that the proof of the exact estimate for $m=3$  turns out to be non-trivial. For $m \geqslant 4$, we were unable to prove that the estimates are exact.

\theorem \label{theorem2} The following equalities hold:

\begin{equation}d_1(T^2) = d_2(T^2) = \frac{\sqrt{2}}{2} = 0.707107...\end{equation}

\begin{equation} d_3(T^2) = \frac{\sqrt{13}}{6} = 0.600925... \label{tor3} \end{equation}

We also highlight several estimates obtained by considering square grids on the torus and searching for colorings in the form of a Boolean satisfiability problem (SAT). In some cases, it is possible to improve the result by reducing the grid coloring to a simpler geometric construction.

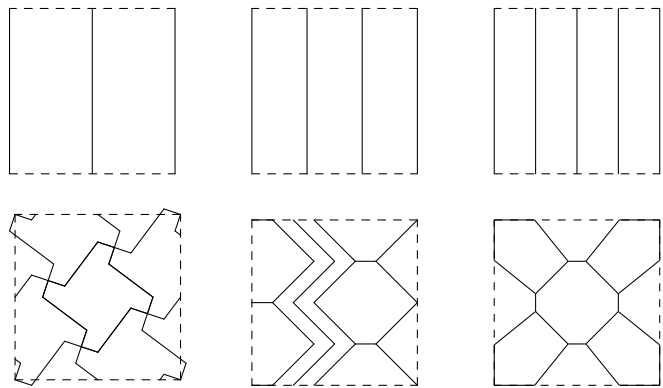
\begin{figure}[!htb]
    \centering
    \begin{tabular}{p{2.8cm} p{2.8cm} p{2.8cm}}
    \definecolor{wwwwww}{rgb}{0.1,0.1,0.1}
\begin{tikzpicture}[line cap=round,line join=round,x=1.0cm,y=1.0cm, scale = 0.22]
\draw[dashed] (0, 0)-- (10, 0);
\draw (10, 0)-- (10, 10);
\draw[dashed] (10, 10)-- (0, 10);
\draw (0, 10)-- (0, 0);
\draw (5, 0) -- (5, 10);
\end{tikzpicture} & \definecolor{wwwwww}{rgb}{0.1,0.1,0.1}
\begin{tikzpicture}[line cap=round,line join=round,x=1.0cm,y=1.0cm, scale = 0.22]
\draw[dashed] (0, 0) -- (10, 0);
\draw (10, 0) -- (10, 10);
\draw[dashed] (10, 10) -- (0, 10);
\draw (0, 10) -- (0, 0);
\draw (3.333, 0) -- (3.333, 10);
\draw (6.666, 0) -- (6.666, 10);
\end{tikzpicture} & \definecolor{wwwwww}{rgb}{0.1,0.1,0.1}
\begin{tikzpicture}[line cap=round,line join=round,x=1.0cm,y=1.0cm, scale = 0.22]
\draw[dashed] (0, 0) -- (10, 0);
\draw (10, 0) -- (10, 10);
\draw[dashed] (10, 10) -- (0, 10);
\draw (0, 10) -- (0, 0);
\draw (2.5, 0) -- (2.5, 10);
\draw (5.0, 0) -- (5.0, 10);
\draw (7.5, 0) -- (7.5, 10);
\end{tikzpicture} \\[0.3cm]
    \definecolor{wwwwww}{rgb}{0.1,0.1,0.1}

\begin{tikzpicture}[line cap=round,line join=round,x=1.0cm,y=1.0cm, scale=2.2]

\draw[dashed] (0, 0) -- (1, 0);
\draw[dashed] (1, 0) -- (1, 1);
\draw[dashed] (1, 1) -- (0, 1);
\draw[dashed] (0, 1) -- (0, 0);

    \draw (0.400000, 0.200000) -- (0.500000, 0.166667);
    \draw (0.500000, 0.166667) -- (0.700000, 0.433333);
    \draw (0.700000, 0.433333) -- (0.800000, 0.400000);
    \draw (0.800000, 0.400000) -- (0.766667, 0.300000);
    \draw (0.766667, 0.300000) -- (1.033333, 0.100000);
    \draw (1.033333, 0.100000) -- (1.000000, 0.000000);
    \draw (1.000000, 0.000000) -- (0.900000, 0.033333);
    \draw (0.900000, 0.033333) -- (0.875000, 0.000000);
    \draw (0.500000, -0.000000) -- (0.366667, 0.100000);
    \draw (0.366667, 0.100000) -- (0.400000, 0.200000);
    \draw (-0.000000, 1.000000) -- (0.100000, 0.966667);
    \draw (0.100000, 0.966667) -- (0.125, 1.000); 
    \draw (0.500000, 1.000000) -- (0.633333, 0.900000);
    \draw (0.633333, 0.900000) -- (0.600000, 0.800000);
    \draw (0.600000, 0.800000) -- (0.500000, 0.833333);
    \draw (0.500000, 0.833333) -- (0.300000, 0.566667);
    \draw (0.300000, 0.566667) -- (0.200000, 0.600000);
    \draw (0.200000, 0.600000) -- (0.233333, 0.700000);
    \draw (0.233333, 0.700000) -- (-0.033333, 0.900000);
    \draw (-0.033333, 0.900000) -- (-0.000000, 1.000000);
    \draw (0.200000, 0.600000) -- (0.300000, 0.566667);
    \draw (0.300000, 0.566667) -- (0.500000, 0.833333);
    \draw (0.500000, 0.833333) -- (0.600000, 0.800000);
    \draw (0.600000, 0.800000) -- (0.566667, 0.700000);
    \draw (0.566667, 0.700000) -- (0.833333, 0.500000);
    \draw (0.833333, 0.500000) -- (0.800000, 0.400000);
    \draw (0.800000, 0.400000) -- (0.700000, 0.433333);
    \draw (0.700000, 0.433333) -- (0.500000, 0.166667);
    \draw (0.500000, 0.166667) -- (0.400000, 0.200000);
    \draw (0.400000, 0.200000) -- (0.433333, 0.300000);
    \draw (0.433333, 0.300000) -- (0.166667, 0.500000);
    \draw (0.166667, 0.500000) -- (0.200000, 0.600000);
    \draw (-0.000000, 0.500000) -- (0.100000, 0.633333);
    \draw (0.100000, 0.633333) -- (0.200000, 0.600000);
    \draw (0.200000, 0.600000) -- (0.166667, 0.500000);
    \draw (0.166667, 0.500000) -- (0.433333, 0.300000);
    \draw (0.433333, 0.300000) -- (0.400000, 0.200000);
    \draw (0.400000, 0.200000) -- (0.300000, 0.233333);
    \draw (0.300000, 0.233333) -- (0.100000, -0.033333);
    \draw (0.100000, -0.033333) -- (0.000000, 0.000000);
    \draw (0.000000, 0.000000) -- (0.033333, 0.100000);
    \draw (0.033333, 0.100000) -- (0, 0.125); 
    \draw (0.600000, 0.800000) -- (0.700000, 0.766667);
    \draw (0.700000, 0.766667) -- (0.900000, 1.033333);
    \draw (0.900000, 1.033333) -- (1.000000, 1.000000);
    \draw (1.000000, 1.000000) -- (0.966667, 0.900000);
    \draw (0.966667, 0.900000) -- (1.000000, 0.875000);
    \draw (1.000000, 0.500000) -- (0.900000, 0.366667);
    \draw (0.900000, 0.366667) -- (0.800000, 0.400000);
    \draw (0.800000, 0.400000) -- (0.833333, 0.500000);
    \draw (0.833333, 0.500000) -- (0.566667, 0.700000);
    \draw (0.566667, 0.700000) -- (0.600000, 0.800000);
\end{tikzpicture} & \definecolor{wwwwww}{rgb}{0.1,0.1,0.1}
\begin{tikzpicture}[line cap=round,line join=round,x=1.0cm,y=1.0cm, scale = 0.275]
\draw[dashed] (0, 0) -- (8, 0);
\draw[dashed] (8, 0) -- (8, 8);
\draw[dashed] (8, 8) -- (0, 8);
\draw[dashed] (0, 8) -- (0, 0);
\draw (0, 4) -- (1, 4);
\draw (1, 0) -- (3, 2);
\draw (3, 2) -- (1, 4);
\draw (1, 4) -- (3, 6);
\draw (3, 6) -- (1, 8);
\draw (2, 0) -- (4, 2);
\draw (4, 2) -- (2, 4);
\draw (2, 4) -- (4, 6);
\draw (4, 6) -- (2, 8);
\draw (3, 0) -- (5, 2);
\draw (5, 2) -- (3, 4);
\draw (3, 4) -- (5, 6);
\draw (5, 6) -- (3, 8);
\draw (5, 2) -- (6, 2);
\draw (5, 6) -- (6, 6);
\draw (8, 0) -- (6, 2);
\draw (6, 2) -- (8, 4);
\draw (8, 4) -- (6, 6);
\draw (6, 6) -- (8, 8);
\draw (0, 0) -- (1, 0);
\draw (0, 8) -- (1, 8);
\end{tikzpicture} & \definecolor{wwwwww}{rgb}{0.1,0.1,0.1}
\begin{tikzpicture}[x=1.0cm,y=1.0cm,scale=2.2]
\draw (0.0538126092553825, -0.25) -- (0.25, -0.0538126092553825);
\draw (0.25, -0.0538126092553825) -- (0.25, 0.0538126092553825);
\draw (0.25, 0.0538126092553825) -- (0.0538126092553825, 0.25);
\draw (0.0538126092553825, 0.25) -- (-0.0538126092553825, 0.25);
\draw (-0.0538126092553825, 0.25) -- (-0.25, 0.0538126092553825);
\draw (-0.25, 0.0538126092553825) -- (-0.25, -0.0538126092553825);
\draw (-0.25, -0.0538126092553825) -- (-0.0538126092553825, -0.25);
\draw (-0.0538126092553825, -0.25) -- (0.0538126092553825, -0.25);
\draw [dashed](-0.5, -0.5) -- (-0.5, 0.5);
\draw [dashed](-0.5, 0.5) -- (0.5, 0.5);
\draw [dashed](0.5, 0.5) -- (0.5, -0.5);
\draw [dashed](0.5, -0.5) -- (-0.5, -0.5);
\draw (-0.5, -0.5) -- (-0.5, -0.25572601923713684); 
\draw (-0.5, -0.5) -- (-0.25572601923713684, -0.5); 
\draw (0.5, -0.5) -- (0.5, -0.25572601923713684); 
\draw (0.5, -0.5) -- (0.25572601923713684, -0.5); 
\draw (0.5, 0.5) -- (0.5, 0.25572601923713684); 
\draw (0.5, 0.5) -- (0.25572601923713684, 0.5); 
\draw (-0.5, 0.5) -- (-0.5, 0.25572601923713684); 
\draw (-0.5, 0.5) -- (-0.25572601923713684, 0.5); 
\draw (0.0538126092553825, -0.25) -- (0.25572601923713684, -0.5);
\draw (0.25, -0.0538126092553825) -- (0.5, -0.25572601923713684);
\draw (0.25, 0.0538126092553825) -- (0.5, 0.25572601923713684);
\draw (0.0538126092553825, 0.25) -- (0.25572601923713684, 0.5);
\draw (-0.0538126092553825, 0.25) -- (-0.25572601923713684, 0.5);
\draw (-0.25, 0.0538126092553825) -- (-0.5, 0.25572601923713684);
\draw (-0.25, -0.0538126092553825) -- (-0.5, -0.25572601923713684);
\draw (-0.0538126092553825, -0.25) -- (-0.25572601923713684, -0.5);
\end{tikzpicture} \\
    \end{tabular}
    
    \caption{Best partitions, $m =  2, \dots, 7$ }
    \label{torus_2_7}
\end{figure}

\theorem \label{theorem3} When $4 \leqslant m \leqslant 7$, the following upper and lower estimates hold:

\begin{equation} 0.556799... = \frac{\sqrt {12401}}{200}  \leqslant d_4(T^2) \leqslant \frac{\sqrt{5}}{4} = 0.559016... \label{tor4} \end{equation}

\begin{equation} 0.521178... = \frac{\sqrt {7850}}{170} \leqslant d_5(T^2) \leqslant  
\frac{\sqrt{10}}{6} = 0.527046... \label{tor5} \end{equation}

\begin{equation} 0.474667... = \frac{\sqrt {73}}{18} \leqslant d_6(T^2) \leqslant \frac{\sqrt{17}}{8} = 0.515388... \label{tor6} \end{equation}

\begin{equation} 0.444444... = \frac{4}{9} \leqslant d_7(T^2) \leqslant \frac{5 - \sqrt{7}}{3} = 0.511452... \label{tor7} \end{equation}

The upper estimates for $2 \leqslant m \leqslant 7$ are provided by the partitions shown in Fig. \ref{torus_2_7}. Table \ref{table:results} summarizes the estimates of $d_m(T^2)$ in $1 \leqslant m \leqslant 25$. The last column contains the ``gap'' between the estimates, which is calculated as the ratio of the difference between the upper and lower estimates and the lower estimate. The table also includes estimates obtained by solving optimization problems, which will be discussed in the following sections.

\setlength{\tabcolsep}{12pt}
\renewcommand{\arraystretch}{0.99}
\begin{center}
\begin{table}[ht]
\caption{Estimates of $d_m(T^2)$}
\label{table:results}
\centering
 \begin{tabular}{ |c|c|c|c|c|  }
 \hline
 $m$ & lower bound & upper bound & best method & gap \\
 \hline\hline
 1 & 0.707107 & 0.707107 & - & exact \\ 
 \hline
 2  & 0.707107 & 0.707107 & stripes & exact \\
 \hline
 3  & 0.600925 & 0.600925 & stripes & exact \\
 \hline
 4  & 0.556707 & 0.559017 & stripes & 0.004 \\
 \hline
 5  & 0.521178 & 0.527046 & SAT+geometry &  0.0113 \\
 \hline
 6  & 0.474667 &  0.515388 & SAT+geometry & 0.0858 \\
 \hline
 7  & 0.444444 & 0.511452 & SAT+geometry & 0.1508 \\
 \hline
 8  & 0.398942 & 0.441942 & hex. grid & 0.1078 \\
 \hline 
 9 & 0.376126 & 0.417229 & glob. opt. & 0.1093 \\
 \hline
 10 & 0.356825 & 0.402869 & glob. opt. & 0.1290 \\
 \hline
 11 & 0.340219 & 0.388518 & hex. grid & 0.1420 \\
 \hline
 12 & 0.325735 & 0.361111 & hex. grid & 0.1086 \\
 \hline
 13 & 0.312956 & 0.353295 & glob. opt. & 0.1289 \\
 \hline
 14 & 0.301572 & 0.339199 & hex. grid & 0.1248 \\
 \hline
 15 & 0.291346 & 0.320555 & hex. grid & 0.1003 \\
 \hline
 16 & 0.282095 & 0.311729 & glob. opt. & 0.1050 \\
 \hline
 17 & 0.273672 & 0.305490 & glob. opt. & 0.1163 \\
 \hline
 18 & 0.265962 & 0.295031 & glob. opt. & 0.1093\\
 \hline
 19 & 0.258868 & 0.290247 & glob. opt. & 0.1212 \\
 \hline
 20 & 0.252313 & 0.279514 & glob. opt. & 0.1078 \\
 \hline
 21 & 0.246233 & 0.274383 & glob. opt. & 0.1143 \\
 \hline
 22 & 0.240571 & 0.268728 & glob. opt. & 0.1170 \\
 \hline
 23 & 0.235283 & 0.259548 & glob. opt. & 0.1031 \\
 \hline
 24 & 0.230329 & 0.254370 & glob. opt. & 0.1044 \\
 \hline
 25 & 0.225676 & 0.249304 & glob. opt. & 0.1047 \\
 \hline

\end{tabular}
\end{table}
\end{center}

\section{Proofs}

 \subsection{Proof of Theorem \ref{theorem1}}

\noindent(\ref{theorem1}). Draw $m$ vertical lines on the torus $T^2$ at a distance $\frac{1}{m}$ from each other. The resulting partition is the desired partition of the torus into $m$ sets. Each set is a ``strip'' in which any two points along the $x$-axis are at a distance of no more than $\frac{1}{m}$, and along the $y$-axis at a distance of no more than $\frac{1}{2}$. Hence, the diameter of each strip is $\sqrt{\frac{1}{4} + \frac{1}{m^2}}$. 

\noindent(\ref{tor_sqrt}). If the diameter of a set does not exceed $\tau < \frac{1}{2}$, it is bounded by a circle with a diameter of $2\tau < 1$. Thus, we can apply reasoning related to the same problem on the plane. This implies that the area of the set does not exceed the area of a circle with a diameter of $\tau$, i.e., $S \leqslant \pi \frac{\tau^2}{4}$.

Let $S_1, \dots, S_m$ be the areas of sets that cover the torus. We can estimate the total area from below by one, since these sets must cover the torus.

$$ S(T^2) = 1 \leqslant S_1 + ... + S_m \leqslant m \cdot \frac{\pi \tau^2}{4}.$$

Hence, $d_m(T^2) \geqslant \tau \geqslant \frac{2}{\sqrt{\pi m}}$ as required.

 \noindent(\ref{tor_lower_lines}). Draw $k$ vertical lines on the flat torus so that the distance between neighboring lines is $\frac{1}{k}$. Suppose that the covering sets have diameter $\tau < \frac{1}{k}$. Consider each of the vertical lines. By the pigeonhole principle, on each of the lines there must be at least $k+1$ different sets from the covering. Since the distance between the sets is at least $\frac{1}{k}$, all these sets on different lines will also be different. Thus, to cover the torus with sets of diameter not exceeding $\tau < \frac{1}{k}$ at least $k(k+1)$ sets are required. It follows that $d_{k^2+k- 1}(T^2) \geqslant \frac{1}{k}$.
 \qed

\subsection{Proof of Theorem \ref{theorem2}}

In this section, as before, we assume that all sums of coordinates in the torus are computed modulo $1$.

Upper bounds for of $m = 2, 3$ follow from estimate (\ref{tor_stripes}) with stripe covering in Theorem \ref{theorem1}. The lower bound for $m=2$ is trivial. It is sufficient to prove the lower estimate for $m=3$.

Let $\tau_3 = \sqrt{\frac{1}{4} + \frac{1}{9}} = \frac{\sqrt{13}}{6} = 0.600925...$. The estimate $d_3(T^2) \leqslant \tau_3$ follows from the covering with three stripes (\ref{tor_stripes}). Let us prove that this value cannot be improved. Consider an arbitrary covering of the torus $T^2$ with three closed sets $F_1, F_2, F_3$ as a coloring of $T^2$ in three colors $1, 2, 3$ respectively. Here we assume that each point has one to three colors. 


\definition The set $C(u) = \{i \in \{1, 2, 3\} \;|\; u \in F_i\}$ is a multicolor of the point $u = (x, y) \in T^2$ (i.e.  indices of sets covering the point \(u\)).


\definition Denote by  \[l^V_x = \{(x, y) \, | \, 0 \leqslant y \leqslant 1\}\] a vertical line, and similarly by \[l^H_y = \{(x, y)\, | \, 0 \leqslant x \leqslant 1\}\] a horizontal line.

\definition Denote by $C(A) = \bigcup\limits_{u \in A} C(u)$ the set of colors  that occur in the set $A$.


\definition Let us call the line $l^V_x$ trichromatic if $C(l^V_x) = \{1, 2, 3\}$, bichromatic $\{p,q\}$-line if $C(l^V_x) = \{p, q\}$, and monochromatic $\{q\}$-line if $C(l^V_x) = \{q\}$


\definition Denote a vertical or a horizontal strip by 
\[l^V_{[\alpha, \beta]} = \{(x, y) \, | \, \alpha \leqslant x \leqslant \beta\},  \quad l^H_{[\alpha, \beta]} = \{(x, y)\, | \, \alpha \leqslant y \leqslant \beta\}.\]

\definition 
Let $\mu_i(l^V_x) = \mu(F_i \cap l^V_x)$, where \(\mu\) is the one-dimensional Lebesgue measure. Note that these sets are measurable, since we are only considering closed sets.

Here we formulate the classical result of Raikov \cite{raikov1939addition} which was later generalized by Macbeath \cite{macbeath1953measure}. In the mentioned articles these statements are formulated for measurable sets, but for our purposes, the case of closed sets is sufficient.
\begin{theorem}[Raikov, \cite{raikov1939addition}]    
Suppose that $A,B$ are closed subsets of $T^1$. Then
    \[\mu(A + B) \geqslant \min(1, \mu(A) + \mu(B)),\]

where $A+B = \{x+y \; | \; x \in A, y \in B\}$ is the Minkowski sum.
    \label{macbeath}
\end{theorem}

\statement \label{trapezoid_lemma} Suppose that \[\exists x: \quad C(l^V_x)=C(l^V_{x + 1/2})=\{1,2,3\}.\] Then there exists a covering set whose diameter is at least $\tau_3$.

\begin{figure}[!htb]
    \centering
    \begin{tabular}{cc}
    \includegraphics[scale=0.14]{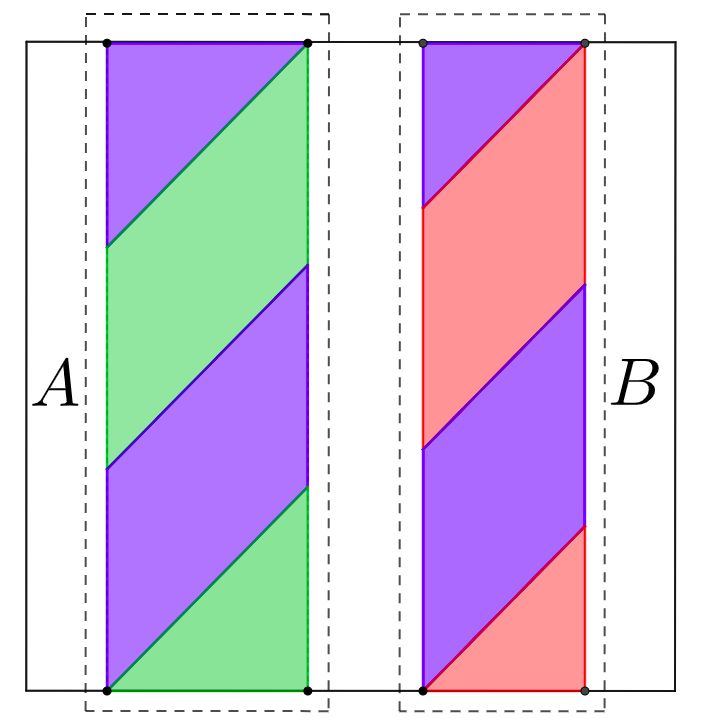} &
    \includegraphics[scale=0.138]{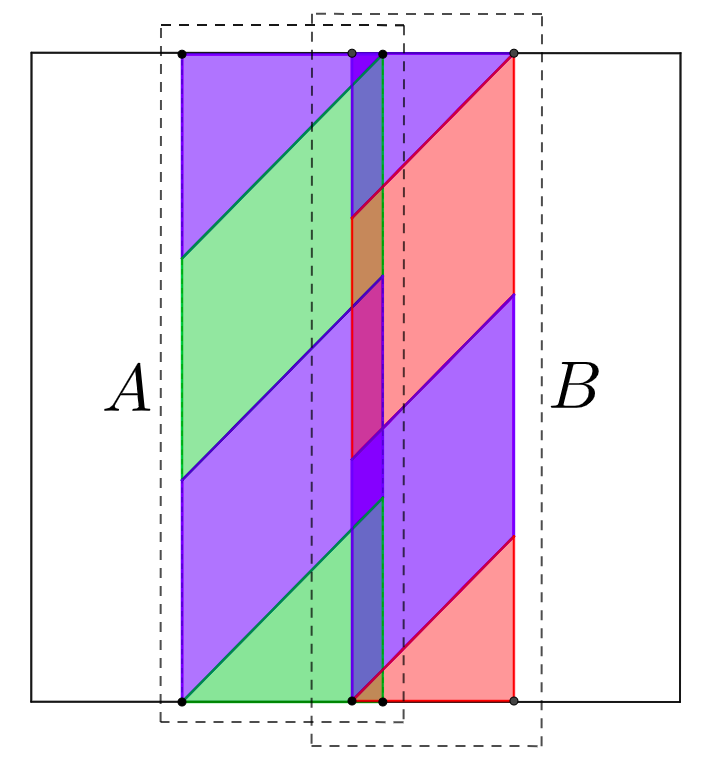} 
    \end{tabular}
    \caption{Arrangement of the green-blue strip (A) and the blue-red strip (B).  1 ~--- red, 2 ~--- blue, 3 ~--- green.}
    \label{d3_strips_location}
\end{figure}

\proof
Without loss of generality, let us prove this for lines $l^V_0$, $l^V_{1/2}$ (Fig. \ref{d3_strips_location}, left). Let $$T^2 \subseteq F_1 \cup F_2 \cup F_3,$$ and 
\[
    \max \{\operatorname{diam} (F_1), \operatorname{diam} (F_2), \operatorname{diam} (F_3) \} = \tau'_3<\tau_3.
\]

Let \(\tau'_3 = \sqrt{\frac{1}{4} + (\frac{1}{3} - \varepsilon)^2}\) for some \(\varepsilon > 0\). We have $$2 = \mu(l^V_0) + \mu(l^V_{1/2}) \leqslant \mu_1(l^V_0) + \mu_2(l^V_0) + \mu_3(l^V_0) + \mu_1(l^V_{1/2}) + \mu_2(l^V_{1/2}) + \mu_3(l^V_{1/2}).$$

Therefore,  
\begin{equation}
\exists i:\; \mu_i(l^V_0) + \mu_i(l^V_{1/2}) \geqslant \frac{2}{3}.
\label{meas12}
\end{equation}

On the other hand, let  $A_i = l^V_0 \cap F_i, B_i = l^V_{1/2} \backslash F_i$. 
For $(x,y) \in A_i$ we have the following:  

$$\left\{ \left(x + \frac{1}{2}, y + \beta\right) \;\Big|\; \frac{1}{3} - \varepsilon \leqslant \beta \leqslant \frac{2}{3} + \varepsilon \right\} \subset B_i.$$ 
  Therefore,  
$$A_i + \left\{\left(\frac{1}{2}, \beta \right) \; \Big| \; \frac{1}{3} - \varepsilon \leqslant \beta \leqslant \frac{2}{3} + \varepsilon\right\} \subset B_i.$$ 
Using Theorem \ref{macbeath}, we can assert 
$$\mu(B_i) \geqslant \min\left(1, \mu(A_i) + \frac{1}{3} + 2\varepsilon\right).$$  Hence, we obtain 
$$\mu(B_i) = 1 - \mu_i(l^V_{1/2}) \geqslant \mu(A_i) + \frac{1}{3} + 2\varepsilon,$$
$$\frac{2}{3} > \mu(A_i) + \mu_i(l^V_{1/2}) = \mu_i(l^V_0) + \mu_i(l^V_{1/2}),$$ a contradiction with \eqref{meas12}.
\qed

\statement \label{trapezoid_lemma_twosets} Let there exist two bichromatic $\{p,q\}$-lines $l^V_x$ and $l^V_{x + \gamma}$, $\frac{1}{3} \leqslant \gamma \leqslant \frac{1}{2}$. Then there exists a covering set whose diameter is at least $\tau_3$.


\proof Consider the covering of two lines $l^V_x, l^V_{x + \gamma}$ as a bichromatic coloring (Fig. \ref{d3_strips_location}, right). Due to the closedness of the covering sets, there exists a bichromatic point $(x, y)$ on the line $l^V_x$ (there are colors $a$ and $b$). At the same time, the point $(x + \gamma, y + \frac{1}{2})$ on the line $l^V_{x + \gamma}$ has color $a$ or $b$ (possibly both). Therefore, the pair of points $(x, y)$ and $(x + \gamma, y + \frac{1}{2})$ is covered by a single set, which means that its diameter is at least $\sqrt{\gamma^2 + \frac{1}{4}} \geqslant \tau_3$. 
\qed

\statement Let there exists a monochromatic $\{c\}$-line $l^V_x$, $c \in \{1,2,3\}$. Then there exists a covering set whose diameter is at least $\tau_3$.

\proof Without loss of generality, let us prove this for the $\{1\}$-line $l^V_0$. In this case, any point on the line $l^V_{1/3}$ or $l^V_{2/3}$ has a color different from $1$.  So $l^V_{1/3}$ and $l^V_{2/3}$, are bichromatic or monochromatic lines. Let us consider the cases of the chromaticities of these lines. Let one of these lines be monochromatic and the other be bichromatic. Without loss of generality, let $l^V_{1/3}$ be a monochromatic $\{2\}$-line and $l^V_{2/3}$ be a bichromatic $\{2,3\}$-line. We choose an arbitrary point $(\frac{2}{3}, y)$ on the line $l^V_{2/3}$ of color $2$. Its distance to the point $(\frac{1}{3}, y + \frac{1}{2})$ is exactly $\tau_3$, and its ends are of the same color, that is, they lie in the same set, and we get the required. If both lines are bichromatic, then by Proposition (\ref{trapezoid_lemma_twosets}) we get the desired result. Finally, if both lines are monochromatic, then they must all have different colors. But then the lines $l^H_0$ and $l^H_{1/2}$ are trichromatic lines.  By Proposition (\ref{trapezoid_lemma}) we get the required. \qed

\begin{figure}[!htb]
    \centering

    \begin{tabular}{ccc}
    \includegraphics[scale=0.26]{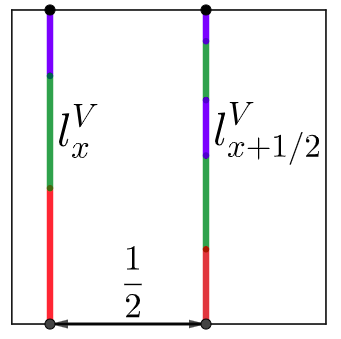} &
    \includegraphics[scale=0.325]{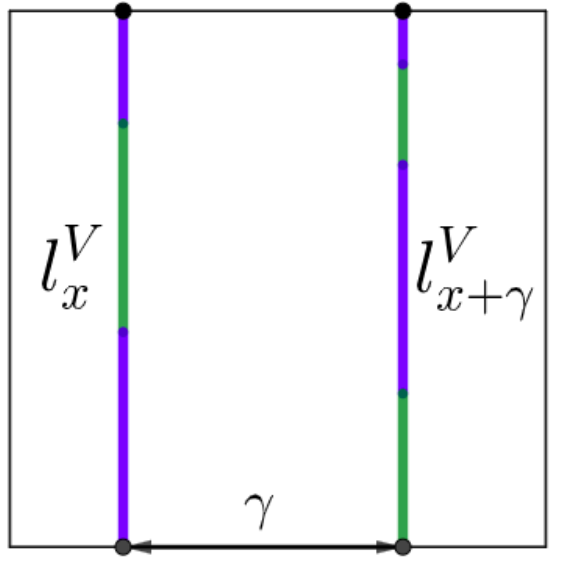} 
    \end{tabular}
    \caption{Propositions 3 and 4}
    \label{conditions_for_the_existence}
\end{figure}

\statement Let the torus $T^2$ be covered by three sets of diameters smaller than $\tau_3$. Fix two different colors $p, q$. Then all bichromatic $\{p,q\}$-lines are enclosed in a strip of width $\frac{1}{3}$. That is, there exists a number $\eta$ such that every bichromatic $\{p,q\}$-line $l^V_x$ has $x \in [\eta, \eta + \frac{1}{3}]$. 

\proof Let us choose an arbitrary $\{p,q\}$-line $l^V_x$. By Proposition 2, every other $\{p,q\}$-line $v_{x^{\prime}}$ has coordinates $x^{\prime} \in [x - \frac{1}{3}, x + \frac{1}{3}]$. Among all $\{p,q\}$-lines with coordinates in the segment $[x - \frac{1}{3}, x + \frac{1}{3}]$ we choose the lines with minimum and maximum coordinates $x_l, x_r$, respectively. By Proposition 2, we have $x_r - x_l \leqslant \frac{1}{3}$. Then all bichromatic $\{p,q\}$-lines lie in the strip $l^V_{[x_l, x_r]}$ having width at most $\frac{1}{3}$. \qed 

\remark It is sufficient to prove the theorem for all mesh colorings (that is, the colorings considered in Section 3.2.2). Then the proof for an arbitrary coloring will follow because of the transition to the limit.

\statement \label {one_color} Let the torus $T^2$ be covered by three sets of diameters smaller than $\tau_3$.  Then some color is present on the whole axis, i.e., there exists $c \in \{1, 2, 3\}$ for which either $\forall x \exists y: c \in C(x, y)$ or $\forall y \exists x: c \in C(x, y)$.

\proof Consider a coloring where the grid lines are divided into $s$ equal segments on each axis, i.e., the  torus $T^2$ is divided into $s^2$ equal square ``pixels''. Pick any arbitrary pixel, let it be of color $c_1$, and consider its connectivity component. Suppose that the color $c_1$ is not present on the entire axis. Then its connectivity component is bounded on both axes. Let this component contain $a_1$ pixels. Then there must be the same color $c_2$ around the boundary of this connected component, otherwise a trichromatic point is found. Consider a connected component of color $c_2$, let its area be $a_2$. This will continue until $a_1 + a_2 + ... + a_m$ does not exceed $s^2$,  a contradiction. \qed

\statement $d_3(T^2) \geqslant \tau_3$ 

\proof  Let the torus $T^2$ be covered by three sets of diameters smaller than $\tau_3$. By Proposition (\ref{one_color}) we see that there exists a color \(c\) present on all lines \(l^V_x, 0 \leqslant x \leqslant 1\). Without loss of generality, let \(c = 2\), then every vertical line \(l^V_x\) has color $\{1,2\}$, $\{2,3\}$, or $\{1,2,3\}$. Let all $\{2,3\}$-lines be in the strip $A = l^V_{[a, a + 1/3]}$, all $\{1,2\}$-lines be in the strip $B = l^V_{[b, b + 1/3]}$. Without loss of generality, let $a < b$.  If $A$ and $B$ do not intersect (Fig. \ref{d3_strips_location}), that is, $a + \frac{1}{3} < b$, then for $t = \frac{1}{2}(a + \frac{1}{3} + b)$ the lines $l^V_m$ and $l^V_{t + 1/2}$ are trichromatic, which contradicts Proposition (\ref{trapezoid_lemma}). If the strips $A$ and $B$ intersect (Fig. \ref{d3_strips_location}), then there exists a trichromatic line $l^V_x \in A \cap B$ (we consider closed sets in the covering). Thus the line $l^V_{x + 1/2} \notin A \cup B$, so it is trichromatic. By Proposition (\ref{trapezoid_lemma}) we get a contradiction.  \qed

  \subsection{Verification of the estimates in Theorem \ref{theorem3}}

\subsubsection{Upper bounds}

Upper bound for $m =  4$ follow from estimate (\ref{tor_stripes}) with stripe covering in Theorem \ref{theorem1}. Upper bounds of $d_5(T^2), d_6(T^2), d_7(T^2)$ were proven by the analysis of the covering from our GitHub repo (see \cite{our_github}). The validation script checks that the considerable covering consists of regions with the specified diameter. 


    

In the case of a torus region bounded by a polyline, it is not enough to consider pairwise distances between vertices.  Distances from each vertex to the boundary for which some coordinate differs by $1/2$ are also required.

\begin{lemma}
    Let
    \[
        d_0(F) = \max_{i\neq j} |P_i P_j|,
    \]
    \[
        Q_i = F \cap \left(\{(x,y) \; | \;  x = x_i \pm 1/2 \} \cup \partial F \cap \{(x,y) \; | \; y = y_i \pm 1/2 \}\right)
    \]
    \[
        d_i(F) = \max_{Q \in Q_i } |P_i Q|, \quad 1 \leqslant i \leqslant N.
    \]
    Then the diameter of $F$ can be computed as follows:
    \[
        \operatorname{diam} F = \max_{0 \leqslant i \leqslant N} d_i(F).
    \]

    \label{diam_lemma_1_over_2}
\end{lemma}

It suffices to show that a figure bounded by a polyline always has a diameter, one end of which is a vertex. Consider the points $u=(x_u, y_u)$, $v = (x_v,y_v)$. Without loss of generality, suppose that the shortest segment $uv$ does not intersect the lines $x=0$, $y=0$, and  $|y_u-y_v|<\frac{1}{2}$.

\medskip

\noindent \emph{Case 1}. Suppose that $v$  lies in the interior of $F$. Then, for a sufficiently small $\varepsilon$ we have $v'= (x_v, y_v+\varepsilon(y_v-y_u)) \in F$, $|uv'|>|uv|$. 

\medskip

\noindent  \emph{Case 2}. $u \in P_i P_{i+1}$, $v \in P_j P_{j+1}$, $|x_u-x_v|<\frac{1}{2}$. If $\varepsilon$ is small enough, then set $v' = v+\varepsilon(P_i - P_j) \in F$, $v'' = v-\varepsilon(P_i - P_j) \in F$. Observe that either $|uv'|>|uv|$ or $|uv''|>|uv|$.

\medskip

\noindent  \emph{Case 3}. $u \in P_i P_{i+1}$, $v \in P_j P_{j+1}$, $|x_u-x_v| = \frac{1}{2}$.  Consider the points $u' = P_i P_{i+1} \cap \{(x,y): x = u_x-\varepsilon\}$, $v' = P_j P_{j+1} \cap \{(x,y): x = v_x-\varepsilon\}$, and $u''= u - uu'$, $v'' = v - vv'$. As before, suppose that $\varepsilon$ is small enough and $u',u'', v', v'' \in F$.  If $uu'=vv'$ then it is possible to shift the segment $uv$ so that one of its ends coincides with the vertex of the polyline. Finally,  if $uu' \neq vv'$, then either $|u'v'|>|uv|$ or $|u''v''|>|uv|$, since in one of the cases the difference of $y$ coordinates  is increasing.
\qed

\begin{figure}[!htb]
    \centering
    \begin{tabular}{p{2.8cm}}
    \definecolor{wwwwww}{rgb}{0.1,0.1,0.1}

\begin{tikzpicture}[line cap=round,line join=round,x=1.0cm,y=1.0cm, scale=3.6]

\draw[dashed] (0.5, -0.2) -- (0.5, 0.67);
\draw[dashed] (-0.15, 0.5) -- (0.65, 0.5);

\draw (0.5, -0.02) node[anchor=north west, font=\footnotesize] {$\{(x_i+1/2,y)\}$};
\draw (0.6, 0.55) node[anchor=north west, font=\footnotesize] {$\{(x, y_i+1/2)\}$};
\draw (-0.25, -0.02) node[anchor=north west, font=\footnotesize] {$P_i = (x_i, y_i)$};

\draw (0,0)--(0.06,0.27);
\draw (0.06,0.27)--(-0.07,0.37);
\draw (-0.07,0.37)--(0.1,0.55);
\draw (0.1,0.55)--(0.2,0.43);
\draw (0.2,0.43)--(0.33,0.54);
\draw (0.33,0.54)--(0.58,0.24);
\draw (0.58,0.24)--(0.2,0.18);
\draw (0.2,0.18)--(0,0);

\fill (0,0) circle (0.4 pt);
\fill (0.0527, 0.5) circle (0.4 pt);
\fill (0.1416, 0.5) circle (0.4 pt);
\fill (0.282, 0.5) circle (0.4 pt);
\fill (0.358, 0.5) circle (0.4 pt);
\fill (0.5, 0.336) circle (0.4 pt);
\fill (0.5, 0.227) circle (0.4 pt);

\end{tikzpicture}\\
    \end{tabular}
    \caption{Definition of the set $Q_i$ in Lemma 1}
    \label{lemma_polygon}
\end{figure}
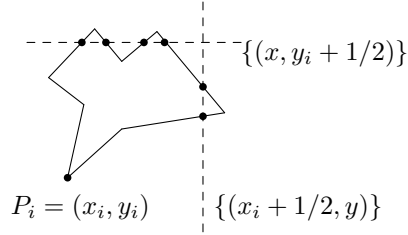

\subsubsection{Lower and upper estimates via SAT approach}

In this section, we consider the coloring of a mesh on a torus in which two square cells must have different colors if the diameter of their union is greater than $\tau$. The corresponding graph is close to the well-known family of Borsuk graphs, which are defined on the sphere $S^{n-1} \subset \mathbb{R}^n$ \cite{Lovasz1983,kahle2020chromatic}.

Let the integer $s \geqslant 2$ be the grid size, and let $\tau > 0$ be the diameter constraint. Denote the vertices and edges of the graph $G_{s, \tau} = \left(V_s, E_{s, \tau} \right)$ in the following way:

$$V_m = \left\{\left(\frac{x}{s}, \frac{y}{s}\right) : x, y \in \{0, ..., s - 1\}  \right\} \hspace{0.5cm} E_{s, \tau} = \left\{ (v_1, v_2) \in V_s \times V_s : \rho_T(v_1, v_2) \geqslant \tau \right\}$$

\statement If the graph $G_{s, \tau}$ can be properly colored with $m$ colors, then
$d_m(T^2) \leqslant \tau $.

\proof If we color the square $\left[\frac{x - 1/2}{s}, \frac{x + 1/2}{s}\right] \times \left[\frac{y - 1/2}{s}, \frac{y + 1/2}{s}\right] \subset T^2$ in the same color as the vertex $\left(\frac{x}{s}, \frac{y}{s} \right) \in V_s$ in the graph coloring, then we obtain the torus partitioning into $m$ parts with diameters not greater than $\tau$.
\qed



\statement If $G_{s, \tau}$ cannot be properly colored in $m$ colors, then $d_m(T^2)\geqslant \tau - \frac{\sqrt{2}}{s}$.

\proof

Consider an arbitrary partition of the torus into $m$ parts and color each of the points of the set $V_m$ in a color equal to the part number (or one of the numbers) to which it belongs. Since the vertices of the graph $G_{s, \tau}$ cannot be properly colored, there exist two points $v_1, v_2 \in V_m$ of the same color, such that $(v_1, v_2) \in E_{s, \tau}$. Then $\rho_T(v_1, v_2)\geqslant\tau$ and these points lie in one part of the partition. This means that the diameter of this part is not less than $\tau$. Without loss of generality, we get that in any torus partition into $m$ parts there exist the part whose diameter is no less than $\tau$. Hence, if $G_{s, \tau}$ cannot be properly colored in $m$ colors, then $d_m(T^2)\geqslant \tau$. \qed 

These facts allow us to compute the lower and upper bounds of the $d_m(T^2)$ using the SAT solver. Here, as in a number of other works, we use the most natural way of transforming the graph coloring problem into a Boolean formula \cite{heule_str}. Namely, for $m$ colors and $s^2$ vertices we have $ms^2$ Boolean variables $w_{uc}$, each of them has the meaning that ``vertex $u$ is colored in color $c$''.  For each edge $(u,v)\in E$ and each color, we write the condition $\neg w_{uc} \lor \neg w_{vc}$. Each vertex must have at least one color, i.e. $w_{u1} \lor w_{u2} \lor \dots \lor w_{um}$. Finally, all the above formulas are joined by a conjunction. As a result, we get a conjunctive normal form (CNF), which is the input to the SAT solver.




\begin{figure}[!htb]
    \centering
    \begin{tabular}{p{3.5cm} p{3.5cm}}
    {\includegraphics[width=0.25\textwidth]{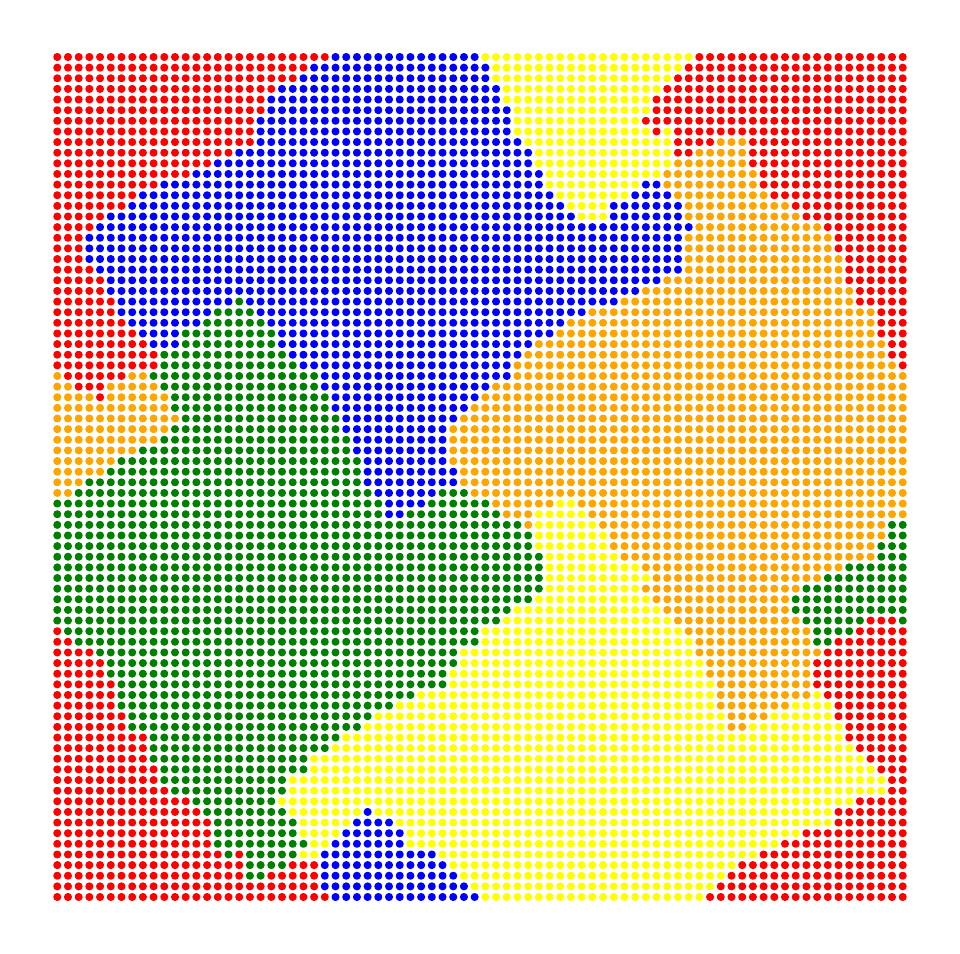}} &
    {\includegraphics[width=0.25\textwidth]{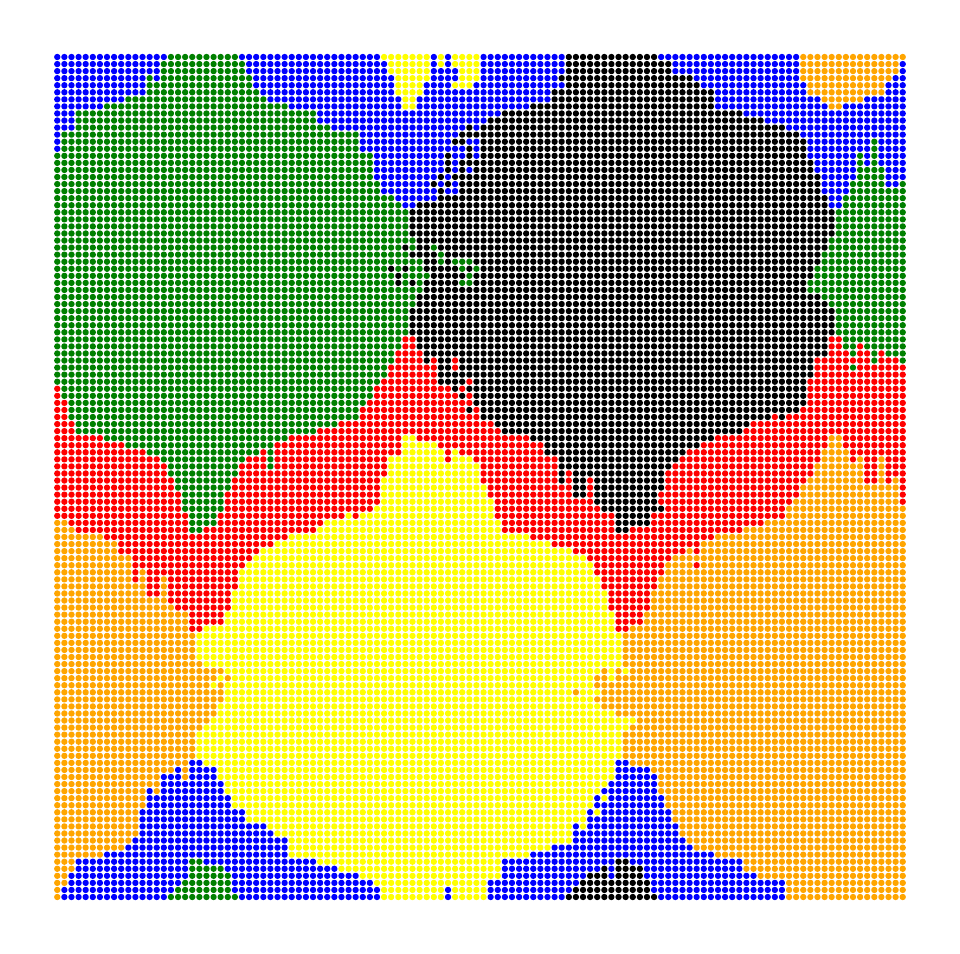}} \\
    \end{tabular}
    \caption{Colorings in 5 and 6 colors found via the \emph{kissat} solver}
    \label{SAT_coloring}
\end{figure}


Table \ref{table:graph_parameters} summarizes the parameters of graphs $G_{s, \tau}$, which cannot be properly colored in $m$ colors. The variable $k$ is the square of the maximum distance at which the SAT solver outputs ``UNSAT'', i.e. the graph cannot be properly colored. The last two rows of Table \ref{table:graph_parameters} correspond to the cases where we did not improve the lower bound of $d_6(T^2)$ due to the hardness of the corresponding SAT problem.

Additionally, the structure of  the proper graph colorings  found by the SAT solver (Fig. \ref{SAT_coloring}) allows us to construct  partitions into 4, 5, 6, 7 parts which brings us to new upper bounds in these cases (Fig. \ref{torus_2_7}).

The \emph{kissat} solver, various modifications of which show the best performance  in SAT competitions, was used to find colorings or to prove their non-existence \cite{Biere2}.

Note that an approach to the Borsuk problem based on discretization of sets was proposed in \cite{zong2}. Moreover, similar methods were used in paper \cite{heule_str} for the coloring of the strip in the setting of the Hadwiger--Nelson problem.  

\setlength{\tabcolsep}{9pt}
\renewcommand{\arraystretch}{0.99}
\begin{center}
\begin{table}[ht]
\caption{Parameters of graphs that can't be properly colored}
\label{table:graph_parameters}
\centering
 \begin{tabular}{ |c|c|c|c|c|c|c|c|c|c|c|c|c|c|  }
 \hline
 $m$ & $s$ (grid size) & $k$ & $\tau = \frac{\sqrt{k}}{s}$ & Time \\
 \hline\hline
 4 & 200 & 12401 & 0.5568 & $\approx 1 h.$ \\
 \hline
 5 & 170 & 7850 & 0.5212 & $\approx 2 h.$ \\
 \hline
 6 & 18 & 73 & 0.4747 & $\approx 5 h.$ \\
 \hline
 7 & 9 & 16 & 0.4444 & $\approx 100 h.$  \\
 \hline
 \hline
 6 & 24 & 137 & 0.4877 &  $>500 h.$ \\
 \hline
 6 & 30 & 221 & 0.4955 & $>500h.$ \\
 \hline

\end{tabular}
\end{table}
\end{center}

\subsubsection{Upper estimates via periodic hexagonal grids}

\begin{figure}[!htb]
    \centering
    \begin{tabular}{p{3cm} p{3cm}  p{3cm}} 
    \definecolor{wwwwww}{rgb}{0.1,0.1,0.1}
\begin{tikzpicture}[x=1.0cm,y=1.0cm,scale=2.4]
\draw (-0.0625, 0.1875) -- (0.0, 0.0);

 \draw[->, line width=0.8pt] (0.0, 0.0) -- (0.1875, -0.0625) node[midway, below=-0.06cm] {\tiny $p_1$};
  \draw[->, line width=0.8pt] (0.1875, -0.0625) -- (0.375, 0.125) node[midway, below right=-0.36cm] {\tiny $p_2$};
  \draw[->, line width=0.8pt] (0.375, 0.125) -- (0.3125, 0.3125) node[midway, right=-0.06cm] {\tiny $p_3$};

\draw (0.3125, 0.3125) -- (0.125, 0.375);
\draw (0.125, 0.375) -- (-0.0625, 0.1875);
\draw (0.0625, 0.5625) -- (0.125, 0.375);
\draw (0.125, 0.375) -- (0.3125, 0.3125);
\draw (0.3125, 0.3125) -- (0.5, 0.5);
\draw (0.5, 0.5) -- (0.4375, 0.6875);
\draw (0.4375, 0.6875) -- (0.25, 0.75);
\draw (0.25, 0.75) -- (0.0625, 0.5625);
\draw (0.1875, 0.9375) -- (0.25, 0.75);
\draw (0.25, 0.75) -- (0.4375, 0.6875);
\draw (0.4375, 0.6875) -- (0.625, 0.875);
\draw (0.625, 0.875) -- (0.5625, 1.0625);
\draw (0.5625, 1.0625) -- (0.375, 1.125);
\draw (0.375, 1.125) -- (0.1875, 0.9375);
\draw (0.3125, 0.3125) -- (0.375, 0.125);
\draw (0.375, 0.125) -- (0.5625, 0.0625);
\draw (0.5625, 0.0625) -- (0.75, 0.25);
\draw (0.75, 0.25) -- (0.6875, 0.4375);
\draw (0.6875, 0.4375) -- (0.5, 0.5);
\draw (0.5, 0.5) -- (0.3125, 0.3125);
\draw (0.4375, 0.6875) -- (0.5, 0.5);
\draw (0.5, 0.5) -- (0.6875, 0.4375);
\draw (0.6875, 0.4375) -- (0.875, 0.625);
\draw (0.875, 0.625) -- (0.8125, 0.8125);
\draw (0.8125, 0.8125) -- (0.625, 0.875);
\draw (0.625, 0.875) -- (0.4375, 0.6875);
\draw (0.5625, 1.0625) -- (0.625, 0.875);
\draw (0.625, 0.875) -- (0.8125, 0.8125);
\draw (0.8125, 0.8125) -- (1.0, 1.0);
\draw (1.0, 1.0) -- (0.9375, 1.1875);
\draw (0.9375, 1.1875) -- (0.75, 1.25);
\draw (0.75, 1.25) -- (0.5625, 1.0625);
\draw (0.6875, 0.4375) -- (0.75, 0.25);
\draw (0.75, 0.25) -- (0.9375, 0.1875);
\draw (0.9375, 0.1875) -- (1.125, 0.375);
\draw (1.125, 0.375) -- (1.0625, 0.5625);
\draw (1.0625, 0.5625) -- (0.875, 0.625);
\draw (0.875, 0.625) -- (0.6875, 0.4375);
\draw (0.8125, 0.8125) -- (0.875, 0.625);
\draw (0.875, 0.625) -- (1.0625, 0.5625);
\draw (1.0625, 0.5625) -- (1.25, 0.75);
\draw (1.25, 0.75) -- (1.1875, 0.9375);
\draw (1.1875, 0.9375) -- (1.0, 1.0);
\draw (1.0, 1.0) -- (0.8125, 0.8125);
\draw[dashed] (0, 0) -- (1, 0);
\draw[dashed] (1, 0) -- (1, 1);
\draw[dashed] (1, 1) -- (0, 1);
\draw[dashed] (0, 1) -- (0, 0);
\end{tikzpicture} &
    \definecolor{wwwwww}{rgb}{0.1,0.1,0.1}
\begin{tikzpicture}[x=1.0cm,y=1.0cm,scale=2.4]
\draw (0.2479338842373698, -0.08264462788311436) -- (0.0, 0.0);
\draw (0.0, 0.0) -- (-0.02479338848990289, 0.09917355393506747);
\draw (-0.02479338848990289, 0.09917355393506747) -- (0.09090909090909091, 0.2727272727272727);
\draw (0.09090909090909091, 0.2727272727272727) -- (0.3388429751464607, 0.19008264484415835);
\draw (0.3388429751464607, 0.19008264484415835) -- (0.36363636363636365, 0.09090909090909088);
\draw (0.36363636363636365, 0.09090909090909088) -- (0.2479338842373698, -0.08264462788311436);
\draw (0.6115702478737335, 0.008264463025976526) -- (0.36363636363636365, 0.09090909090909088);
\draw (0.36363636363636365, 0.09090909090909088) -- (0.3388429751464608, 0.19008264484415835);
\draw (0.3388429751464608, 0.19008264484415835) -- (0.4545454545454546, 0.3636363636363636);
\draw (0.4545454545454546, 0.3636363636363636) -- (0.7024793387828243, 0.2809917357532492);
\draw (0.7024793387828243, 0.2809917357532492) -- (0.7272727272727273, 0.18181818181818177);
\draw (0.7272727272727273, 0.18181818181818177) -- (0.6115702478737335, 0.008264463025976526);
\draw (0.9752066115100971, 0.09917355393506741) -- (0.7272727272727273, 0.18181818181818177);
\draw (0.7272727272727273, 0.18181818181818177) -- (0.7024793387828244, 0.2809917357532492);
\draw (0.7024793387828244, 0.2809917357532492) -- (0.8181818181818182, 0.4545454545454545);
\draw (0.8181818181818182, 0.4545454545454545) -- (1.066115702419188, 0.37190082666234014);
\draw (1.066115702419188, 0.37190082666234014) -- (1.0909090909090908, 0.27272727272727265);
\draw (1.0909090909090908, 0.27272727272727265) -- (0.9752066115100971, 0.09917355393506741);
\draw (0.3388429751464608, 0.19008264484415835) -- (0.09090909090909094, 0.2727272727272727);
\draw (0.09090909090909094, 0.2727272727272727) -- (0.06611570241918804, 0.37190082666234014);
\draw (0.06611570241918804, 0.37190082666234014) -- (0.18181818181818185, 0.5454545454545454);
\draw (0.18181818181818185, 0.5454545454545454) -- (0.42975206605555166, 0.4628099175714311);
\draw (0.42975206605555166, 0.4628099175714311) -- (0.4545454545454546, 0.3636363636363636);
\draw (0.4545454545454546, 0.3636363636363636) -- (0.3388429751464608, 0.19008264484415835);
\draw (0.7024793387828244, 0.2809917357532492) -- (0.4545454545454546, 0.3636363636363636);
\draw (0.4545454545454546, 0.3636363636363636) -- (0.4297520660555517, 0.4628099175714311);
\draw (0.4297520660555517, 0.4628099175714311) -- (0.5454545454545455, 0.6363636363636362);
\draw (0.5454545454545455, 0.6363636363636362) -- (0.7933884296919154, 0.5537190084805219);
\draw (0.7933884296919154, 0.5537190084805219) -- (0.8181818181818182, 0.4545454545454545);
\draw (0.8181818181818182, 0.4545454545454545) -- (0.7024793387828244, 0.2809917357532492);
\draw (1.0661157024191878, 0.37190082666234014) -- (0.8181818181818181, 0.4545454545454545);
\draw (0.8181818181818181, 0.4545454545454545) -- (0.7933884296919153, 0.5537190084805219);
\draw (0.7933884296919153, 0.5537190084805219) -- (0.9090909090909091, 0.7272727272727272);
\draw (0.9090909090909091, 0.7272727272727272) -- (1.157024793328279, 0.6446280993896129);
\draw (1.157024793328279, 0.6446280993896129) -- (1.1818181818181817, 0.5454545454545454);
\draw (1.1818181818181817, 0.5454545454545454) -- (1.0661157024191878, 0.37190082666234014);
\draw (0.4297520660555517, 0.4628099175714311) -- (0.18181818181818188, 0.5454545454545454);
\draw (0.18181818181818188, 0.5454545454545454) -- (0.15702479332827898, 0.6446280993896129);
\draw (0.15702479332827898, 0.6446280993896129) -- (0.2727272727272728, 0.8181818181818181);
\draw (0.2727272727272728, 0.8181818181818181) -- (0.5206611569646427, 0.7355371902987038);
\draw (0.5206611569646427, 0.7355371902987038) -- (0.5454545454545455, 0.6363636363636362);
\draw (0.5454545454545455, 0.6363636363636362) -- (0.4297520660555517, 0.4628099175714311);
\draw (0.7933884296919153, 0.5537190084805219) -- (0.5454545454545454, 0.6363636363636362);
\draw (0.5454545454545454, 0.6363636363636362) -- (0.5206611569646425, 0.7355371902987037);
\draw (0.5206611569646425, 0.7355371902987037) -- (0.6363636363636364, 0.909090909090909);
\draw (0.6363636363636364, 0.909090909090909) -- (0.8842975206010062, 0.8264462812077946);
\draw (0.8842975206010062, 0.8264462812077946) -- (0.9090909090909091, 0.7272727272727271);
\draw (0.9090909090909091, 0.7272727272727271) -- (0.7933884296919153, 0.5537190084805219);
\draw (1.157024793328279, 0.6446280993896129) -- (0.9090909090909092, 0.7272727272727272);
\draw (0.9090909090909092, 0.7272727272727272) -- (0.8842975206010063, 0.8264462812077946);
\draw (0.8842975206010063, 0.8264462812077946) -- (1.0, 0.9999999999999999);
\draw (1.0, 0.9999999999999999) -- (1.24793388423737, 0.9173553721168856);
\draw (1.24793388423737, 0.9173553721168856) -- (1.272727272727273, 0.8181818181818181);
\draw (1.272727272727273, 0.8181818181818181) -- (1.157024793328279, 0.6446280993896129);
\draw (0.5206611569646425, 0.7355371902987038) -- (0.2727272727272727, 0.8181818181818181);
\draw (0.2727272727272727, 0.8181818181818181) -- (0.2479338842373698, 0.9173553721168856);
\draw (0.2479338842373698, 0.9173553721168856) -- (0.36363636363636365, 1.0909090909090908);
\draw (0.36363636363636365, 1.0909090909090908) -- (0.6115702478737335, 1.0082644630259765);
\draw (0.6115702478737335, 1.0082644630259765) -- (0.6363636363636364, 0.909090909090909);
\draw (0.6363636363636364, 0.909090909090909) -- (0.5206611569646425, 0.7355371902987038);
\draw (0.8842975206010063, 0.8264462812077946) -- (0.6363636363636365, 0.909090909090909);
\draw (0.6363636363636365, 0.909090909090909) -- (0.6115702478737336, 1.0082644630259765);
\draw (0.6115702478737336, 1.0082644630259765) -- (0.7272727272727274, 1.1818181818181817);
\draw (0.7272727272727274, 1.1818181818181817) -- (0.9752066115100972, 1.0991735539350673);
\draw (0.9752066115100972, 1.0991735539350673) -- (1.0, 0.9999999999999998);
\draw (1.0, 0.9999999999999998) -- (0.8842975206010063, 0.8264462812077946);
\draw[dashed] (0, 0) -- (1, 0);
\draw[dashed] (1, 0) -- (1, 1);
\draw[dashed] (1, 1) -- (0, 1);
\draw[dashed] (0, 1) -- (0, 0);
\end{tikzpicture} &
    \definecolor{wwwwww}{rgb}{0.1,0.1,0.1}
\begin{tikzpicture}[x=1.0cm,y=1.0cm,scale=2.4]
\draw (0.11111111111911347, -0.16666666668855512) -- (0.0, 0.0);
\draw (0.0, 0.0) -- (0.11111111111911347, 0.1666666666447782);
\draw (0.11111111111911347, 0.1666666666447782) -- (0.25, 0.16666666666666666);
\draw (0.25, 0.16666666666666666) -- (0.3611111111191135, -2.1888463264119196e-11);
\draw (0.3611111111191135, -2.1888463264119196e-11) -- (0.25, -0.16666666666666666);
\draw (0.25, -0.16666666666666666) -- (0.11111111111911347, -0.16666666668855512);
\draw (0.11111111111911347, 0.1666666666447782) -- (0.0, 0.3333333333333333);
\draw (0.0, 0.3333333333333333) -- (0.11111111111911347, 0.4999999999781115);
\draw (0.11111111111911347, 0.4999999999781115) -- (0.25, 0.5);
\draw (0.25, 0.5) -- (0.3611111111191135, 0.3333333333114449);
\draw (0.3611111111191135, 0.3333333333114449) -- (0.25, 0.16666666666666666);
\draw (0.25, 0.16666666666666666) -- (0.11111111111911347, 0.1666666666447782);
\draw (0.3611111111191135, -2.1888463264119196e-11) -- (0.25, 0.16666666666666666);
\draw (0.25, 0.16666666666666666) -- (0.3611111111191135, 0.3333333333114449);
\draw (0.3611111111191135, 0.3333333333114449) -- (0.5, 0.3333333333333333);
\draw (0.5, 0.3333333333333333) -- (0.6111111111191134, 0.1666666666447782);
\draw (0.6111111111191134, 0.1666666666447782) -- (0.5, 0.0);
\draw (0.5, 0.0) -- (0.3611111111191135, -2.1888463264119196e-11);
\draw (0.6111111111191134, -0.16666666668855512) -- (0.5, 0.0);
\draw (0.5, 0.0) -- (0.6111111111191134, 0.1666666666447782);
\draw (0.6111111111191134, 0.1666666666447782) -- (0.75, 0.16666666666666666);
\draw (0.75, 0.16666666666666666) -- (0.8611111111191134, -2.1888463264119196e-11);
\draw (0.8611111111191134, -2.1888463264119196e-11) -- (0.75, -0.16666666666666666);
\draw (0.75, -0.16666666666666666) -- (0.6111111111191134, -0.16666666668855512);
\draw (0.11111111111911347, 0.4999999999781115) -- (0.0, 0.6666666666666666);
\draw (0.0, 0.6666666666666666) -- (0.11111111111911347, 0.8333333333114448);
\draw (0.11111111111911347, 0.8333333333114448) -- (0.25, 0.8333333333333333);
\draw (0.25, 0.8333333333333333) -- (0.3611111111191135, 0.6666666666447781);
\draw (0.3611111111191135, 0.6666666666447781) -- (0.25, 0.5);
\draw (0.25, 0.5) -- (0.11111111111911347, 0.4999999999781115);
\draw (0.3611111111191135, 0.3333333333114449) -- (0.25, 0.5);
\draw (0.25, 0.5) -- (0.3611111111191135, 0.6666666666447782);
\draw (0.3611111111191135, 0.6666666666447782) -- (0.5, 0.6666666666666666);
\draw (0.5, 0.6666666666666666) -- (0.6111111111191134, 0.4999999999781115);
\draw (0.6111111111191134, 0.4999999999781115) -- (0.5, 0.33333333333333337);
\draw (0.5, 0.33333333333333337) -- (0.3611111111191135, 0.3333333333114449);
\draw (0.6111111111191134, 0.1666666666447782) -- (0.5, 0.3333333333333333);
\draw (0.5, 0.3333333333333333) -- (0.6111111111191134, 0.4999999999781115);
\draw (0.6111111111191134, 0.4999999999781115) -- (0.75, 0.5);
\draw (0.75, 0.5) -- (0.8611111111191134, 0.3333333333114449);
\draw (0.8611111111191134, 0.3333333333114449) -- (0.75, 0.16666666666666666);
\draw (0.75, 0.16666666666666666) -- (0.6111111111191134, 0.1666666666447782);
\draw (0.8611111111191134, -2.188849101969481e-11) -- (0.75, 0.16666666666666663);
\draw (0.75, 0.16666666666666663) -- (0.8611111111191134, 0.3333333333114448);
\draw (0.8611111111191134, 0.3333333333114448) -- (1.0, 0.33333333333333326);
\draw (1.0, 0.33333333333333326) -- (1.1111111111191134, 0.16666666664477817);
\draw (1.1111111111191134, 0.16666666664477817) -- (1.0, -2.7755575615628914e-17);
\draw (1.0, -2.7755575615628914e-17) -- (0.8611111111191134, -2.188849101969481e-11);
\draw (0.3611111111191135, 0.6666666666447782) -- (0.25, 0.8333333333333334);
\draw (0.25, 0.8333333333333334) -- (0.3611111111191135, 0.9999999999781115);
\draw (0.3611111111191135, 0.9999999999781115) -- (0.5, 1.0);
\draw (0.5, 1.0) -- (0.6111111111191134, 0.8333333333114449);
\draw (0.6111111111191134, 0.8333333333114449) -- (0.5, 0.6666666666666667);
\draw (0.5, 0.6666666666666667) -- (0.3611111111191135, 0.6666666666447782);
\draw (0.6111111111191134, 0.4999999999781116) -- (0.5, 0.6666666666666667);
\draw (0.5, 0.6666666666666667) -- (0.6111111111191134, 0.833333333311445);
\draw (0.6111111111191134, 0.833333333311445) -- (0.75, 0.8333333333333334);
\draw (0.75, 0.8333333333333334) -- (0.8611111111191134, 0.6666666666447782);
\draw (0.8611111111191134, 0.6666666666447782) -- (0.75, 0.5000000000000001);
\draw (0.75, 0.5000000000000001) -- (0.6111111111191134, 0.4999999999781116);
\draw (0.8611111111191134, 0.3333333333114449) -- (0.75, 0.5);
\draw (0.75, 0.5) -- (0.8611111111191134, 0.6666666666447782);
\draw (0.8611111111191134, 0.6666666666447782) -- (1.0, 0.6666666666666666);
\draw (1.0, 0.6666666666666666) -- (1.1111111111191134, 0.4999999999781115);
\draw (1.1111111111191134, 0.4999999999781115) -- (1.0, 0.33333333333333337);
\draw (1.0, 0.33333333333333337) -- (0.8611111111191134, 0.3333333333114449);
\draw (0.8611111111191134, 0.6666666666447781) -- (0.75, 0.8333333333333333);
\draw (0.75, 0.8333333333333333) -- (0.8611111111191134, 0.9999999999781115);
\draw (0.8611111111191134, 0.9999999999781115) -- (1.0, 0.9999999999999999);
\draw (1.0, 0.9999999999999999) -- (1.1111111111191134, 0.8333333333114448);
\draw (1.1111111111191134, 0.8333333333114448) -- (1.0, 0.6666666666666666);
\draw (1.0, 0.6666666666666666) -- (0.8611111111191134, 0.6666666666447781);
\draw[dashed] (0, 0) -- (1, 0);
\draw[dashed] (1, 0) -- (1, 1);
\draw[dashed] (1, 1) -- (0, 1);
\draw[dashed] (0, 1) -- (0, 0);
\end{tikzpicture} \\
    \end{tabular}
    \caption{Hexagonal tilings, $m=8, 11, 12$}
    \label{torus_lattice}
\end{figure}
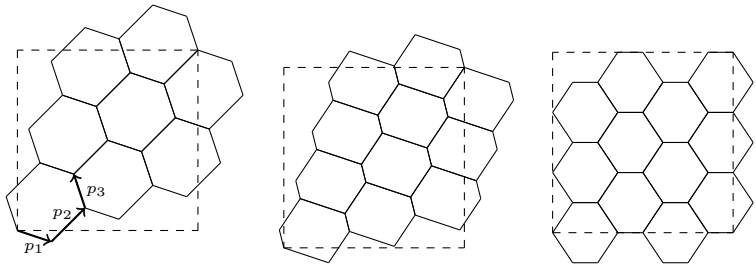

In this section, we consider the hexagonal tilings on the flat torus (Fig. \ref{torus_lattice}). Suppose that the sets $F_i$ are congruent to a centrally symmetric hexagon $H$. Let $\vec{p}_1 = (x_1, y_1)$, $\vec{p}_2 = (x_2, y_2)$, $\vec{p}_3 = (x_3, y_3)$ be vectors corresponding to three consecutive sides of $H$. In the case of partition into $m \geqslant 7$ parts, the diameter of $H$ is the length of its longest diagonal. 

Let $\vec{h}_1 = (1, 0)$, $\vec{h}_2 = (0, 1)$ be vectors corresponding to the sides of the square (flat torus). Observe that $\vec{h}_1$ and $\vec{h}_2$ can be expressed as follows: 

\begin{equation}
\centering
\left\{\begin{aligned}
a_1 \vec{p}_1 + a_2 \vec{p}_2 + a_3 \vec{p}_3 = \vec{h}_1, \\
b_1 \vec{p}_1 + b_2 \vec{p}_2 + b_3 \vec{p}_3 = \vec{h}_2. \\ 
\end{aligned}\right.
\label{grid_system}
\end{equation}




The general solution of \ref{grid_system} can be written in the following form by expressing variables $x_2, x_3$ by $x_1$ and variables $y_2, y_3$ by $y_1$.

\begin{equation} 
    \begin{cases}
    x_2 = -x_1 + \alpha_2 \\
    x_3 = x_1 + \alpha_3 \\
    y_2 = -y_1 + \beta_2 \\
    y_3 = y_1 + \beta_3, \\
    \end{cases}
    \label {x23y23}
\end{equation}

where 

$$\alpha_2 = \frac{b_3}{a_2 b_3 - a_3 b_2}, \alpha_3 = \frac{-b_2}{a_2 b_3 - a_3 b_2}, \beta_2 = \frac{-a_3}{a_2 b_3 - a_3 b_2}, \beta_3 = \frac{a_2}{a_2 b_3 - a_3 b_2}.$$



Consider the system \ref{grid_system} and assume that $\vec{p}_3 = (0, 0)$. We then compute the determinant for the left and right parts of \ref{grid_system} we get the following:
$$(a_1b_2 - a_2b_1) \cdot \Delta(\vec{p}_1, \vec{p}_2) = \Delta(\vec{h}_1, \vec{h}_2) = 1,$$
where $\Delta(u,v)$ is the determinant of the $2 \times 2$ matrix with columns $u,v$.

If $\vec{p}_3 = (0, 0)$, then the flat torus is divided into parallelograms based on the vectors $\vec{p}_1, \vec{p}_2$. Since the area of such parallelogram is $\Delta(\vec{p}_1, \vec{p}_2)$, we obtain  that $|a_1b_2 - a_2b_1| = m$, where $m$ is the number of parts. We also proved that $\sgn (a_1b_2 - a_2b_1) = \sgn \Delta(\vec{p}_1, \vec{p}_2)$. Finally, $$a_1b_2 - a_2b_1 = m \cdot \sgn \Delta(\vec{p}_1, \vec{p}_2).$$

Similarly, it follows that 
$$a_2 b_3 - a_3 b_2 = m \cdot \Delta(\vec{p}_2, \vec{p}_3), \quad a_3b_1 - a_1b_3 = m \cdot \sgn \Delta(\vec{p}_3, \vec{p}_1).$$ Since $\vec{p}_1, \vec{p}_2, \vec{p}_3$ correspond to three consecutive sides of this hexagon $H$, we have $$ \sgn \Delta(\vec{p}_1, \vec{p}_2) = \sgn \Delta(\vec{p}_1, \vec{p}_2) = -\sgn \Delta(\vec{p}_3, \vec{p}_1). $$ 
Thus, the equations $$a_1 b_2 - a_2 b_1 = a_2 b_3 - a_3 b_2 = a_1 b_3-a_3 b_1   = \pm m$$
are satisfied for any hexagonal tiling of the torus.

Let us denote by $d(x_1, y_1)$ the diameter of the hexagon $H$, where the coordinates $\vec{p}_1 = (x_1, y_1)$, $\vec{p}_2$, and $\vec{p}_3$ are calculated by \eqref{x23y23}.

Since the diameter of a hexagon $H$ is the length of its longest diagonal, $d(x_1, y_1)$ can be rewritten in the following form:

$$d(x_1, y_1) = \max\left(|\vec{p}_1 + \vec{p}_2 + \vec{p}_3|, |\vec{p}_2 + \vec{p}_3 - \vec{p}_1|, |\vec{p}_3 - \vec{p}_1 - \vec{p}_2|\right).$$




Denote by $f: \mathbb{R}^2 \rightarrow \mathbb{R}_+$ the function such that $f(x_1, y_1) = d^2(x_1, y_1)$. 

\begin {equation}
f(x_1, y_1) =  x_1^2 + y_1^2 + C + 2 \cdot \max 
\begin{pmatrix}
 |x_1(\alpha_2 + \alpha_3) + y_1(\beta_2 + \beta_3)| + (\alpha_2 \alpha_3 + \beta_2 \beta_3), \\
 x_1(\alpha_3 - \alpha_2) + y_1(\beta_3 - \beta_2) - (\alpha_2 \alpha_3 + \beta_2 \beta_3)
\end{pmatrix},
\end {equation}
where $C = \alpha_2^2 + \alpha_3^2 + \beta_2^2 + \beta_3^2$ is the constant.

Let \(\tau = ({ \displaystyle\min_{x_1, y_1} f(x_1, y_1)})^{1/2}\), and \((x^*, y^*) = \displaystyle\arg \min_{x_1, y_1} f(x_1, y_1)\). Then, the partition of the torus into $m$ copies of hexagon $H$ based on $\vec{p}_1 = (x^*, y^*)$ brings us to the upper estimate $d_m(T^2) \leqslant \tau$.

Table \ref{table:grid_estimation_results} contains the parameters necessary to obtain upper estimates of the value $d_m(T^2)$ in the case of $7 \leqslant m \leqslant 16$:







\setlength{\tabcolsep}{9pt}
\renewcommand{\arraystretch}{0.99}
\begin{center}
\begin{table}[ht]
\caption{Upper estimates of $d_m(T^2)$}
\label{table:grid_estimation_results}
\centering
 \begin{tabular}{ |c|c|c|c|c|c|c|c|c|c|c|c|c|c|  }
 \hline
 $m$ & $a_1$ & $a_2$ & $a_3$ & $b_1$ & $b_2$ & $b_3$ &  $f_{\min}=\tau^2$ & $x^*$ & $y^*$ \\
 \hline\hline
 7 & 3 & 1 & -2 & -1 & 2 & 3 & 650/2401 & 12/49 & 4/49 \\ 
 \hline
 8 & 3 & 2 & -1 & -1 & 2 & 3 & 25/128 & 3/16 & -1/16 \\
 \hline
 9 & 1 & 3 & 2 & 3 & 0 & -3 & 130/729 & 2/27 & 2/9 \\
 \hline
 10 & 3 & 2 & -1 & -2 & 2 & 4 & 221/1250 & 9/50 & -3/25 \\
 \hline
 11 & -1 & 2 & 3 & 4 & 3 & -1 & 2210/14641 & -3/121 & 12/121 \\
 \hline
 12 & 2 & 4 & 2 & 3 & 0 & -3 & 169/1296 & 1/9 & 1/6 \\
 \hline
 14 & 4 & 2 & -2 & 1 & 4 & 3 & 1105/9604 & 8/49 & 2/49 \\
 \hline
 15 & 4 & 3 & -1 & -1 & 3 & 4 & 578/5625 & 4/25 & -1/25 \\
 \hline
 16 & 4 & 2 & -2 & 0 & 4 & 4 & 25/256 & 3/16 & 0 \\
 \hline

\end{tabular}
\end{table}
\end{center}

\noindent Let us illustrate the above mentioned technique in the example for the upper estimate of $d_8(T^2)$.

\begin{align}
    &\begin{cases}
    3 x_1 + 2 x_2 - x_3 = 1 \\
    3 y_1 + 2 y_2 - y_3 = 0 \\
    - x_1 + 2 x_2 + 3 x_3 = 0 \\
    - y_1 + 2 y_2 + 3 y_3 = 1 \\
    \end{cases} 
    \Longleftrightarrow
    \begin{cases}
    x_2 = -x_1 + \frac{3}{8} \\
    x_3 = x_1 - \frac{2}{8} \\
    y_2 = -y_1 + \frac{1}{8} \\
    y_3 = y_1 + \frac{2}{8} \\
    \end{cases}
\end{align}

Using formulas for $\alpha_2, \alpha_3, \beta_2, \beta_3$ the optimization problem can be written as $$f(x, y) = x^2 + y^2 + \frac{9}{32} + 2 \cdot\max \left(\left|\frac{x}{8}  + \frac{3y}{8} \right| - \frac{1}{16}, \frac{-5x}{8}  + \frac{y}{8}  + \frac{1}{16}\right) \rightarrow \min.$$

Here we omit the computation of minima. The solution is
 \[\tau^2 = \displaystyle \min_{x, y} f(x, y) = f\left(\frac{3}{16},\frac{-1}{16}\right)=\frac{25}{128}.\]
 Thus, there is a partition of the torus into $8$ equal hexagons with diameter $\tau =  \frac{5}{8\sqrt 2}$.  Hence, $d_8(T^2) \leqslant \frac{5}{8\sqrt{2}} = 0.441941...$ (Fig. \ref{torus_lattice}, left).


\subsubsection{Global optimization approach}

Our paper \cite{DAM_our_article} proposes an optimization algorithm for computing the upper estimates of $d_m(F)$ for planar sets, which will be adapted to find the optimal partitions of the flat torus.

\begin{figure}[!htb]
    \centering
    \begin{tabular}{p{3cm} p{3cm}  p{3cm}} 
    \definecolor{wwwwww}{rgb}{0.1,0.1,0.1}
\begin{tikzpicture}[line cap=round,line join=round,x=1.0cm,y=1.0cm, scale = 2.2]
\draw [dashed](0, 0) -- (1, 0);
\draw [dashed](1, 0) -- (1, 1);
\draw [dashed](1, 1) -- (0, 1);
\draw [dashed](0, 1) -- (0, 0);
\draw (0.23,0.31) circle (4.6 pt);
\draw (0.64,0.01) circle (4.6 pt);
\draw (0.06,0.76) circle (4.6 pt);
\draw (0.69,0.68) circle (4.6 pt);
\draw (0.56,0.34) circle (4.6 pt);
\draw (0.36,0.61) circle (4.6 pt);
\draw (0.31,0.98) circle (4.6 pt);
\draw (0.96,0.12) circle (4.6 pt);
\draw (0.93,0.45) circle (4.6 pt);
\fill   (0.23,0.31) circle (0.4 pt);
\fill   (0.64,0.01) circle (0.4 pt);
\fill   (0.06,0.76) circle (0.4 pt);
\fill   (0.69,0.68) circle (0.4 pt);
\fill   (0.56,0.34) circle (0.4 pt);
\fill   (0.36,0.61) circle (0.4 pt);
\fill   (0.31,0.98) circle (0.4 pt);
\fill   (0.96,0.12) circle (0.4 pt);
\fill   (0.93,0.45) circle (0.4 pt);
\draw (0.03, -0.03) node[anchor=north west, font=\small] {$ a) $};
\end{tikzpicture} &
    \definecolor{uququq}{rgb}{0.1,0.1,0.1}
\begin{tikzpicture}[line cap=round,line join=round,x=1.0cm,y=1.0cm, scale = 2.2]
\draw [dashed](0, 0) -- (1, 0);
\draw [dashed](1, 0) -- (1, 1);
\draw [dashed](1, 1) -- (0, 1);
\draw [dashed](0, 1) -- (0, 0);
\draw (0.411,0.179)--(0.389,0.421);
\draw (0.389,0.421)--(0.150,0.529);
\draw (0.150,0.529)--(0.039,0.292);
\draw (0.039,0.292)--(0.161,0.108);
\draw (0.161,0.108)--(0.411,0.179);
\draw (0.830,0.884)--(0.480,0.823);
\draw (0.480,0.823)--(0.474,0.831);
\draw (0.474,0.831)--(0.502,1.102);
\draw (0.502,1.102)--(0.725,1.173);
\draw (0.725,1.173)--(0.855,0.891);
\draw (0.855,0.891)--(0.830,0.884);
\draw (-0.110,0.649)--(0.138,0.580);
\draw (0.138,0.580)--(0.262,0.787);
\draw (0.262,0.787)--(0.098,0.963);
\draw (0.098,0.963)--(-0.145,0.891);
\draw (-0.145,0.891)--(-0.170,0.884);
\draw (-0.170,0.884)--(-0.110,0.649);
\draw (0.480,0.823)--(0.551,0.539);
\draw (0.551,0.539)--(0.722,0.477);
\draw (0.722,0.477)--(0.890,0.649);
\draw (0.890,0.649)--(0.830,0.884);
\draw (0.830,0.884)--(0.480,0.823);
\draw (0.389,0.421)--(0.551,0.539);
\draw (0.551,0.539)--(0.722,0.477);
\draw (0.722,0.477)--(0.813,0.299);
\draw (0.813,0.299)--(0.725,0.173);
\draw (0.725,0.173)--(0.502,0.102);
\draw (0.502,0.102)--(0.411,0.179);
\draw (0.411,0.179)--(0.389,0.421);
\draw (0.480,0.823)--(0.551,0.539);
\draw (0.551,0.539)--(0.389,0.421);
\draw (0.389,0.421)--(0.150,0.529);
\draw (0.150,0.529)--(0.138,0.580);
\draw (0.138,0.580)--(0.262,0.787);
\draw (0.262,0.787)--(0.474,0.831);
\draw (0.474,0.831)--(0.480,0.823);
\draw (0.262,0.787)--(0.474,0.831);
\draw (0.474,0.831)--(0.502,1.102);
\draw (0.502,1.102)--(0.411,1.179);
\draw (0.411,1.179)--(0.161,1.108);
\draw (0.161,1.108)--(0.098,0.963);
\draw (0.098,0.963)--(0.262,0.787);
\draw (0.725,0.173)--(0.813,0.299);
\draw (0.813,0.299)--(1.039,0.292);
\draw (1.039,0.292)--(1.161,0.108);
\draw (1.161,0.108)--(1.098,-0.037);
\draw (1.098,-0.037)--(0.855,-0.109);
\draw (0.855,-0.109)--(0.725,0.173);
\draw (0.890,0.649)--(0.722,0.477);
\draw (0.722,0.477)--(0.813,0.299);
\draw (0.813,0.299)--(1.039,0.292);
\draw (1.039,0.292)--(1.150,0.529);
\draw (1.150,0.529)--(1.138,0.580);
\draw (1.138,0.580)--(0.890,0.649);
\fill   (0.23,0.31) circle (0.4 pt);
\fill   (0.64,1.01) circle (0.4 pt);
\fill   (0.06,0.76) circle (0.4 pt);
\fill   (0.69,0.68) circle (0.4 pt);
\fill   (0.56,0.34) circle (0.4 pt);
\fill   (0.36,0.61) circle (0.4 pt);
\fill   (0.31,0.98) circle (0.4 pt);
\fill   (0.96,0.12) circle (0.4 pt);
\fill   (0.93,0.45) circle (0.4 pt);
\draw (0.03,-0.03) node[anchor=north west, font=\small] {$ b) $};
\end{tikzpicture} &
    \definecolor{uququq}{rgb}{0.1,0.1,0.1}
\begin{tikzpicture}[line cap=round,line join=round,x=1.0cm,y=1.0cm, scale = 2.2]
\draw [dashed](0, 0) -- (1, 0);
\draw [dashed](1, 0) -- (1, 1);
\draw [dashed](1, 1) -- (0, 1);
\draw [dashed](0, 1) -- (0, 0);
\draw (0.430,0.188)--(0.417,0.438);
\draw (0.417,0.438)--(0.176,0.519);
\draw (0.176,0.519)--(0.040,0.279);
\draw (0.040,0.279)--(0.121,0.143);
\draw (0.121,0.143)--(0.430,0.188);
\draw (0.824,0.810)--(0.514,0.764);
\draw (0.514,0.764)--(0.423,0.855);
\draw (0.423,0.855)--(0.522,1.097);
\draw (0.522,1.097)--(0.704,1.136);
\draw (0.704,1.136)--(0.836,0.916);
\draw (0.836,0.916)--(0.824,0.810);
\draw (-0.093,0.673)--(0.125,0.618);
\draw (0.125,0.618)--(0.241,0.816);
\draw (0.241,0.816)--(0.109,1.035);
\draw (0.109,1.035)--(-0.164,0.916);
\draw (-0.164,0.916)--(-0.176,0.810);
\draw (-0.176,0.810)--(-0.093,0.673);
\draw (0.514,0.764)--(0.529,0.514);
\draw (0.529,0.514)--(0.768,0.433);
\draw (0.768,0.433)--(0.907,0.673);
\draw (0.907,0.673)--(0.824,0.810);
\draw (0.824,0.810)--(0.514,0.764);
\draw (0.417,0.438)--(0.529,0.514);
\draw (0.529,0.514)--(0.768,0.433);
\draw (0.768,0.433)--(0.820,0.333);
\draw (0.820,0.333)--(0.704,0.136);
\draw (0.704,0.136)--(0.522,0.097);
\draw (0.522,0.097)--(0.430,0.188);
\draw (0.430,0.188)--(0.417,0.438);
\draw (0.514,0.764)--(0.529,0.514);
\draw (0.529,0.514)--(0.417,0.438);
\draw (0.417,0.438)--(0.176,0.519);
\draw (0.176,0.519)--(0.125,0.618);
\draw (0.125,0.618)--(0.241,0.816);
\draw (0.241,0.816)--(0.423,0.855);
\draw (0.423,0.855)--(0.514,0.764);
\draw (0.522,0.097)--(0.430,0.188);
\draw (0.430,0.188)--(0.121,0.143);
\draw (0.121,0.143)--(0.109,0.035);
\draw (0.109,0.035)--(0.241,-0.184);
\draw (0.241,-0.184)--(0.423,-0.145);
\draw (0.423,-0.145)--(0.522,0.097);
\draw (0.704,0.136)--(0.820,0.333);
\draw (0.820,0.333)--(1.040,0.279);
\draw (1.040,0.279)--(1.121,0.143);
\draw (1.121,0.143)--(1.109,0.035);
\draw (1.109,0.035)--(0.836,-0.084);
\draw (0.836,-0.084)--(0.704,0.136);
\draw (0.907,0.673)--(0.768,0.433);
\draw (0.768,0.433)--(0.820,0.333);
\draw (0.820,0.333)--(1.040,0.279);
\draw (1.040,0.279)--(1.176,0.519);
\draw (1.176,0.519)--(1.125,0.618);
\draw (1.125,0.618)--(0.907,0.673);
\draw (-0.01,-0.03) node[anchor=north west, font=\small] {$ c) $};
\end{tikzpicture}\\
    \end{tabular}
    
    \caption{Approximate circle packing (a), the Voronoi diagram (b), and the final partition (c).}
    \label{alg_example}
\end{figure}
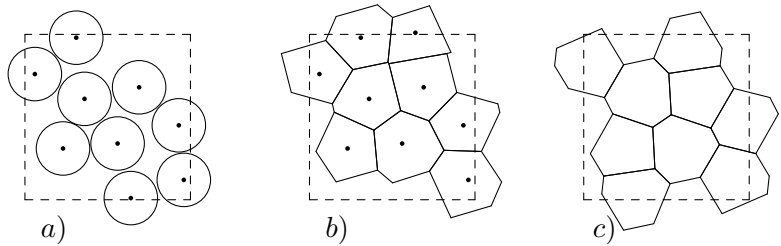

Consider some initial partition of the torus $T^2$ into polygons $T^2 = F_1 \cup F_2 \cup \dots \cup F_k$. Denote vertices of the partition by $X = \{x_1, \dots , x_r \}$. Suppose, in addition, that the number of parts is sufficiently large that the diameters of $F_i$ are less than $1/2$. In this case it is not necessary to consider additional points on the boundary to calculate the diameter, as in Lemma \ref{diam_lemma_1_over_2}.  The vertices of the polygon $F_i$ are the points $x_j$, $j \in \mathcal{I}_i$, where $\mathcal{I}_i \subset \{1,2, \dots, r\}$. We have the following finite-dimensional optimisation problem:
\begin{equation}
\varphi(X) = \max_{i} \operatorname{diam} F_i = \max_{i} \max_{p,q \in \mathcal{I}_i } \| x_p - x_q \| \rightarrow \min.
\label{globopt}
\end{equation}

\begin{figure}[!htb]
    \centering
    \begin{tabular}{p{3cm} p{3cm}  p{3cm}} 
    \definecolor{wwwwww}{rgb}{0.1,0.1,0.1}
\begin{tikzpicture}[x=1.0cm,y=1.0cm,scale=2.4]
\draw (0.2537073386562738, 0.5253057838670099) -- (0.011815163084531729, 0.624039959143586);
\draw (0.29236729340526185, 0.343083424897077) -- (0.2537073386562738, 0.5253057838670099);
\draw (0.07302412838169711, 0.2113239718553244) -- (0.29236729340526185, 0.343083424897077);
\draw (-0.03426938073052299, 0.22339255518827753) -- (0.07302412838169711, 0.2113239718553244);
\draw (-0.0789829328937041, 0.5332969327724387) -- (-0.03426938073052299, 0.22339255518827753);
\draw (0.011815163084531729, 0.624039959143586) -- (-0.0789829328937041, 0.5332969327724387);
\draw (0.29236729340526185, 0.343083424897077) -- (0.07302412838169711, 0.2113239718553244);
\draw (0.4899680900725027, 0.22679446740311987) -- (0.29236729340526185, 0.343083424897077);
\draw (0.437209568449306, 0.012721750847445978) -- (0.4899680900725027, 0.22679446740311987);
\draw (0.29982734523345445, -0.07408130511228639) -- (0.437209568449306, 0.012721750847445978);
\draw (0.19239045398014346, -0.06199287650047549) -- (0.29982734523345445, -0.07408130511228639);
\draw (0.07302412838169711, 0.2113239718553244) -- (0.19239045398014346, -0.06199287650047549);
\draw (0.3448202311066229, 0.6166726134599998) -- (0.29982734523345445, 0.9259186948877136);
\draw (0.2537073386562738, 0.5253057838670099) -- (0.3448202311066229, 0.6166726134599998);
\draw (0.011815163084531729, 0.624039959143586) -- (0.2537073386562738, 0.5253057838670099);
\draw (-0.026867079229921265, 0.8062266422624792) -- (0.011815163084531729, 0.624039959143586);
\draw (0.19239045398014346, 0.9380071234995245) -- (-0.026867079229921265, 0.8062266422624792);
\draw (0.29982734523345445, 0.9259186948877136) -- (0.19239045398014346, 0.9380071234995245);
\draw (0.5901130909426273, 0.27916171857637084) -- (0.4899680900725027, 0.22679446740311987);
\draw (0.8234295303825419, 0.13241308582164169) -- (0.5901130909426273, 0.27916171857637084);
\draw (0.7754465916746461, -0.07748176214025626) -- (0.8234295303825419, 0.13241308582164169);
\draw (0.6757242148310599, -0.1291921646082389) -- (0.7754465916746461, -0.07748176214025626);
\draw (0.437209568449306, 0.012721750847445978) -- (0.6757242148310599, -0.1291921646082389);
\draw (0.4899680900725027, 0.22679446740311987) -- (0.437209568449306, 0.012721750847445978);
\draw (0.8234295303825419, 0.13241308582164169) -- (0.965730619269477, 0.22339255518827753);
\draw (0.7754465916746461, -0.07748176214025626) -- (0.8234295303825419, 0.13241308582164169);
\draw (0.9731329207700787, -0.1937733577375208) -- (0.7754465916746461, -0.07748176214025626);
\draw (1.1923904539801435, -0.06199287650047549) -- (0.9731329207700787, -0.1937733577375208);
\draw (1.073024128381697, 0.2113239718553244) -- (1.1923904539801435, -0.06199287650047549);
\draw (0.965730619269477, 0.22339255518827753) -- (1.073024128381697, 0.2113239718553244);
\draw (0.5946377412236282, 0.6306753117813507) -- (0.3448202311066229, 0.6166726134599998);
\draw (0.6708849440969221, 0.5186012723995379) -- (0.5946377412236282, 0.6306753117813507);
\draw (0.5901130909426273, 0.27916171857637084) -- (0.6708849440969221, 0.5186012723995379);
\draw (0.4899680900725027, 0.22679446740311987) -- (0.5901130909426273, 0.27916171857637084);
\draw (0.29236729340526185, 0.343083424897077) -- (0.4899680900725027, 0.22679446740311987);
\draw (0.2537073386562738, 0.5253057838670099) -- (0.29236729340526185, 0.343083424897077);
\draw (0.3448202311066229, 0.6166726134599998) -- (0.2537073386562738, 0.5253057838670099);
\draw (0.5946377412236282, 0.6306753117813507) -- (0.6757242148310599, 0.8708078353917611);
\draw (0.3448202311066229, 0.6166726134599998) -- (0.5946377412236282, 0.6306753117813507);
\draw (0.29982734523345445, 0.9259186948877136) -- (0.3448202311066229, 0.6166726134599998);
\draw (0.437209568449306, 1.012721750847446) -- (0.29982734523345445, 0.9259186948877136);
\draw (0.6757242148310599, 0.8708078353917611) -- (0.437209568449306, 1.012721750847446);
\draw (0.5946377412236282, 0.6306753117813507) -- (0.6708849440969221, 0.5186012723995379);
\draw (0.6757242148310599, 0.8708078353917611) -- (0.5946377412236282, 0.6306753117813507);
\draw (0.7754465916746461, 0.9225182378597437) -- (0.6757242148310599, 0.8708078353917611);
\draw (0.9731329207700787, 0.8062266422624792) -- (0.7754465916746461, 0.9225182378597437);
\draw (1.0118151630845318, 0.624039959143586) -- (0.9731329207700787, 0.8062266422624792);
\draw (0.9210170671062959, 0.5332969327724387) -- (1.0118151630845318, 0.624039959143586);
\draw (0.6708849440969221, 0.5186012723995379) -- (0.9210170671062959, 0.5332969327724387);
\draw (0.6708849440969221, 0.5186012723995379) -- (0.9210170671062959, 0.5332969327724387);
\draw (0.5901130909426273, 0.27916171857637084) -- (0.6708849440969221, 0.5186012723995379);
\draw (0.8234295303825419, 0.13241308582164169) -- (0.5901130909426273, 0.27916171857637084);
\draw (0.965730619269477, 0.22339255518827753) -- (0.8234295303825419, 0.13241308582164169);
\draw (0.9210170671062959, 0.5332969327724387) -- (0.965730619269477, 0.22339255518827753);
\draw [dashed](0, 0) -- (1, 0);
\draw [dashed](1, 0) -- (1, 1);
\draw [dashed](1, 1) -- (0, 1);
\draw [dashed](0, 1) -- (0, 0);
\draw (-0.02, -0.14) node[anchor=north west] { \scriptsize $d_{ 9 }(T^2) \leqslant 0.417232 $ };
\end{tikzpicture} &
    \definecolor{wwwwww}{rgb}{0.1,0.1,0.1}
\begin{tikzpicture}[x=1.0cm,y=1.0cm,scale=2.4]
\draw (0.7826321977914418, 0.5794505943727307) -- (0.9805534156051805, 0.425658940726029);
\draw (0.7856613071392086, 0.7120058890288267) -- (0.7826321977914418, 0.5794505943727307);
\draw (1.1090417387385925, 0.8075462030389894) -- (0.7856613071392086, 0.7120058890288267);
\draw (1.176314295729081, 0.6652500923169129) -- (1.1090417387385925, 0.8075462030389894);
\draw (1.0530222742996, 0.41056726894304935) -- (1.176314295729081, 0.6652500923169129);
\draw (0.9805534156051805, 0.425658940726029) -- (1.0530222742996, 0.41056726894304935);
\draw (0.23708803137834528, 0.2669365567103549) -- (0.4304123465098062, 0.11481585688962795);
\draw (0.16017648863165784, 0.2825345488164181) -- (0.23708803137834528, 0.2669365567103549);
\draw (0.03723278463377688, 0.026742180390586595) -- (0.16017648863165784, 0.2825345488164181);
\draw (0.10392434783603127, -0.11334564258633917) -- (0.03723278463377688, 0.026742180390586595);
\draw (0.4270430615467923, -0.019341860310478354) -- (0.10392434783603127, -0.11334564258633917);
\draw (0.4304123465098062, 0.11481585688962795) -- (0.4270430615467923, -0.019341860310478354);
\draw (0.44147288296767073, 0.517587285372507) -- (0.23708803137834528, 0.2669365567103549);
\draw (0.39607968481856665, 0.6091792855450616) -- (0.44147288296767073, 0.517587285372507);
\draw (0.1763142957290811, 0.6652500923169129) -- (0.39607968481856665, 0.6091792855450616);
\draw (0.053022274299599964, 0.41056726894304935) -- (0.1763142957290811, 0.6652500923169129);
\draw (0.16017648863165784, 0.2825345488164181) -- (0.053022274299599964, 0.41056726894304935);
\draw (0.23708803137834528, 0.2669365567103549) -- (0.16017648863165784, 0.2825345488164181);
\draw (0.7856613071392086, 0.7120058890288267) -- (0.7063515578609681, 0.8212094631873049);
\draw (1.1090417387385925, 0.8075462030389894) -- (0.7856613071392086, 0.7120058890288267);
\draw (1.1039243478360312, 0.8866543574136608) -- (1.1090417387385925, 0.8075462030389894);
\draw (1.037232784633777, 1.0267421803905865) -- (1.1039243478360312, 0.8866543574136608);
\draw (0.8163633981788718, 1.0844695276160499) -- (1.037232784633777, 1.0267421803905865);
\draw (0.7063515578609681, 0.8212094631873049) -- (0.8163633981788718, 1.0844695276160499);
\draw (0.39607968481856665, 0.6091792855450616) -- (0.1763142957290811, 0.6652500923169129);
\draw (0.5066156398690205, 0.8729842455476036) -- (0.39607968481856665, 0.6091792855450616);
\draw (0.4270430615467923, 0.9806581396895216) -- (0.5066156398690205, 0.8729842455476036);
\draw (0.10392434783603127, 0.8866543574136608) -- (0.4270430615467923, 0.9806581396895216);
\draw (0.10904173873859252, 0.8075462030389894) -- (0.10392434783603127, 0.8866543574136608);
\draw (0.1763142957290811, 0.6652500923169129) -- (0.10904173873859252, 0.8075462030389894);
\draw (0.16017648863165784, 0.2825345488164181) -- (0.053022274299599964, 0.41056726894304935);
\draw (0.03723278463377688, 0.026742180390586595) -- (0.16017648863165784, 0.2825345488164181);
\draw (-0.18363660182112818, 0.08446952761604977) -- (0.03723278463377688, 0.026742180390586595);
\draw (-0.22859436116410714, 0.1766837716007357) -- (-0.18363660182112818, 0.08446952761604977);
\draw (-0.0194465843948195, 0.425658940726029) -- (-0.22859436116410714, 0.1766837716007357);
\draw (0.053022274299599964, 0.41056726894304935) -- (-0.0194465843948195, 0.425658940726029);
\draw (0.44147288296767073, 0.517587285372507) -- (0.5809155563746343, 0.4770073238623984);
\draw (0.23708803137834528, 0.2669365567103549) -- (0.44147288296767073, 0.517587285372507);
\draw (0.4304123465098062, 0.11481585688962795) -- (0.23708803137834528, 0.2669365567103549);
\draw (0.6369366412643866, 0.21772731113298585) -- (0.4304123465098062, 0.11481585688962795);
\draw (0.5809155563746343, 0.4770073238623984) -- (0.6369366412643866, 0.21772731113298585);
\draw (0.5809155563746343, 0.4770073238623984) -- (0.6369366412643866, 0.21772731113298585);
\draw (0.7826321977914418, 0.5794505943727307) -- (0.5809155563746343, 0.4770073238623984);
\draw (0.9805534156051805, 0.425658940726029) -- (0.7826321977914418, 0.5794505943727307);
\draw (0.7714056388358929, 0.1766837716007357) -- (0.9805534156051805, 0.425658940726029);
\draw (0.6369366412643866, 0.21772731113298585) -- (0.7714056388358929, 0.1766837716007357);
\draw (0.6369366412643866, 0.21772731113298585) -- (0.7714056388358929, 0.1766837716007357);
\draw (0.4304123465098062, 0.11481585688962795) -- (0.6369366412643866, 0.21772731113298585);
\draw (0.4270430615467923, -0.019341860310478354) -- (0.4304123465098062, 0.11481585688962795);
\draw (0.5066156398690205, -0.12701575445239643) -- (0.4270430615467923, -0.019341860310478354);
\draw (0.7063515578609681, -0.17879053681269508) -- (0.5066156398690205, -0.12701575445239643);
\draw (0.8163633981788718, 0.08446952761604977) -- (0.7063515578609681, -0.17879053681269508);
\draw (0.7714056388358929, 0.1766837716007357) -- (0.8163633981788718, 0.08446952761604977);
\draw (0.44147288296767073, 0.517587285372507) -- (0.39607968481856665, 0.6091792855450616);
\draw (0.5809155563746343, 0.4770073238623984) -- (0.44147288296767073, 0.517587285372507);
\draw (0.7826321977914418, 0.5794505943727307) -- (0.5809155563746343, 0.4770073238623984);
\draw (0.7856613071392086, 0.7120058890288267) -- (0.7826321977914418, 0.5794505943727307);
\draw (0.7063515578609681, 0.8212094631873049) -- (0.7856613071392086, 0.7120058890288267);
\draw (0.5066156398690205, 0.8729842455476036) -- (0.7063515578609681, 0.8212094631873049);
\draw (0.39607968481856665, 0.6091792855450616) -- (0.5066156398690205, 0.8729842455476036);
\draw[dashed] (0, 0) -- (1, 0);
\draw[dashed] (1, 0) -- (1, 1);
\draw[dashed] (1, 1) -- (0, 1);
\draw[dashed] (0, 1) -- (0, 0);
\draw (-0.02, -0.14) node[anchor=north west] { \scriptsize $d_{ 10 }(T^2) \leqslant 0.402923 $ };
\end{tikzpicture} &
    \definecolor{wwwwww}{rgb}{0.1,0.1,0.1}
\begin{tikzpicture}[x=1.0cm,y=1.0cm,scale=2.4]
\draw (0.3733376304156993, 0.36662118557643525) -- (0.575171729114199, 0.1964736027800764);
\draw (0.3693092562925226, 0.36660844794230335) -- (0.3733376304156993, 0.36662118557643525);
\draw (0.2278857793699955, 0.22551096405712281) -- (0.3693092562925226, 0.36660844794230335);
\draw (0.2741823383090442, 0.03695693264361367) -- (0.2278857793699955, 0.22551096405712281);
\draw (0.3724469383107146, 0.012075238816062317) -- (0.2741823383090442, 0.03695693264361367);
\draw (0.48405748734574877, 0.03624796816824638) -- (0.3724469383107146, 0.012075238816062317);
\draw (0.575171729114199, 0.1964736027800764) -- (0.48405748734574877, 0.03624796816824638);
\draw (0.3693092562925226, 0.36660844794230335) -- (0.19941209400682502, 0.5730758583981066);
\draw (0.2278857793699955, 0.22551096405712281) -- (0.3693092562925226, 0.36660844794230335);
\draw (0.029119958870400562, 0.2682951846009918) -- (0.2278857793699955, 0.22551096405712281);
\draw (0.014785338912755241, 0.3625911016697077) -- (0.029119958870400562, 0.2682951846009918);
\draw (0.0463383542934066, 0.5040581284051607) -- (0.014785338912755241, 0.3625911016697077);
\draw (0.19941209400682502, 0.5730758583981066) -- (0.0463383542934066, 0.5040581284051607);
\draw (0.2278857793699955, 0.22551096405712281) -- (0.029119958870400562, 0.2682951846009918);
\draw (0.2741823383090442, 0.03695693264361367) -- (0.2278857793699955, 0.22551096405712281);
\draw (0.05577406758917024, -0.08137105955691115) -- (0.2741823383090442, 0.03695693264361367);
\draw (-0.0777970055850743, 0.05352141081399571) -- (0.05577406758917024, -0.08137105955691115);
\draw (0.029119958870400562, 0.2682951846009918) -- (-0.0777970055850743, 0.05352141081399571);
\draw (0.19941209400682502, 0.5730758583981066) -- (0.0463383542934066, 0.5040581284051607);
\draw (0.226091137624496, 0.6891457308467227) -- (0.19941209400682502, 0.5730758583981066);
\draw (0.05996880820124079, 0.8546598819906802) -- (0.226091137624496, 0.6891457308467227);
\draw (-0.12819554621557638, 0.7026819692452801) -- (0.05996880820124079, 0.8546598819906802);
\draw (-0.1281152936921236, 0.6736083554769277) -- (-0.12819554621557638, 0.7026819692452801);
\draw (0.0463383542934066, 0.5040581284051607) -- (-0.1281152936921236, 0.6736083554769277);
\draw (0.34352486593015663, 0.7176683875265142) -- (0.3724469383107146, 1.0120752388160623);
\draw (0.226091137624496, 0.6891457308467227) -- (0.34352486593015663, 0.7176683875265142);
\draw (0.05996880820124079, 0.8546598819906802) -- (0.226091137624496, 0.6891457308467227);
\draw (0.05577406758917024, 0.9186289404430888) -- (0.05996880820124079, 0.8546598819906802);
\draw (0.2741823383090442, 1.0369569326436137) -- (0.05577406758917024, 0.9186289404430888);
\draw (0.3724469383107146, 1.0120752388160623) -- (0.2741823383090442, 1.0369569326436137);
\draw (0.5176560761734916, 0.6886283467274534) -- (0.34352486593015663, 0.7176683875265142);
\draw (0.6641883932565097, 0.8689206838898882) -- (0.5176560761734916, 0.6886283467274534);
\draw (0.48405748734574877, 1.0362479681682464) -- (0.6641883932565097, 0.8689206838898882);
\draw (0.3724469383107146, 1.0120752388160623) -- (0.48405748734574877, 1.0362479681682464);
\draw (0.34352486593015663, 0.7176683875265142) -- (0.3724469383107146, 1.0120752388160623);
\draw (0.5520153085255308, 0.5495884614045772) -- (0.5176560761734916, 0.6886283467274534);
\draw (0.3733376304156993, 0.36662118557643525) -- (0.5520153085255308, 0.5495884614045772);
\draw (0.3693092562925226, 0.36660844794230335) -- (0.3733376304156993, 0.36662118557643525);
\draw (0.19941209400682502, 0.5730758583981066) -- (0.3693092562925226, 0.36660844794230335);
\draw (0.226091137624496, 0.6891457308467227) -- (0.19941209400682502, 0.5730758583981066);
\draw (0.34352486593015663, 0.7176683875265142) -- (0.226091137624496, 0.6891457308467227);
\draw (0.5176560761734916, 0.6886283467274534) -- (0.34352486593015663, 0.7176683875265142);
\draw (0.6950107605184371, 0.5157164564996014) -- (0.8718847063078764, 0.6736083554769277);
\draw (0.7274838039916538, 0.3498003862811505) -- (0.6950107605184371, 0.5157164564996014);
\draw (1.0147853389127552, 0.3625911016697077) -- (0.7274838039916538, 0.3498003862811505);
\draw (1.0463383542934066, 0.5040581284051607) -- (1.0147853389127552, 0.3625911016697077);
\draw (0.8718847063078764, 0.6736083554769277) -- (1.0463383542934066, 0.5040581284051607);
\draw (0.575171729114199, 0.1964736027800764) -- (0.48405748734574877, 0.03624796816824638);
\draw (0.6890486869404121, 0.22259505992096848) -- (0.575171729114199, 0.1964736027800764);
\draw (0.8378098994600067, 0.055376375675268395) -- (0.6890486869404121, 0.22259505992096848);
\draw (0.7058167568558498, -0.13111589409341895) -- (0.8378098994600067, 0.055376375675268395);
\draw (0.6641883932565097, -0.1310793161101118) -- (0.7058167568558498, -0.13111589409341895);
\draw (0.48405748734574877, 0.03624796816824638) -- (0.6641883932565097, -0.1310793161101118);
\draw (0.5520153085255308, 0.5495884614045772) -- (0.3733376304156993, 0.36662118557643525);
\draw (0.6950107605184371, 0.5157164564996014) -- (0.5520153085255308, 0.5495884614045772);
\draw (0.7274838039916538, 0.3498003862811505) -- (0.6950107605184371, 0.5157164564996014);
\draw (0.6890486869404121, 0.22259505992096848) -- (0.7274838039916538, 0.3498003862811505);
\draw (0.575171729114199, 0.1964736027800764) -- (0.6890486869404121, 0.22259505992096848);
\draw (0.3733376304156993, 0.36662118557643525) -- (0.575171729114199, 0.1964736027800764);
\draw (0.7274838039916538, 0.3498003862811505) -- (1.0147853389127552, 0.3625911016697077);
\draw (0.6890486869404121, 0.22259505992096848) -- (0.7274838039916538, 0.3498003862811505);
\draw (0.8378098994600067, 0.055376375675268395) -- (0.6890486869404121, 0.22259505992096848);
\draw (0.9222029944149257, 0.05352141081399571) -- (0.8378098994600067, 0.055376375675268395);
\draw (1.0291199588704005, 0.2682951846009918) -- (0.9222029944149257, 0.05352141081399571);
\draw (1.0147853389127552, 0.3625911016697077) -- (1.0291199588704005, 0.2682951846009918);
\draw (0.7058167568558498, 0.868884105906581) -- (0.8378098994600067, 1.0553763756752683);
\draw (0.8718044537844236, 0.7026819692452801) -- (0.7058167568558498, 0.868884105906581);
\draw (1.0599688082012408, 0.8546598819906802) -- (0.8718044537844236, 0.7026819692452801);
\draw (1.0557740675891703, 0.9186289404430888) -- (1.0599688082012408, 0.8546598819906802);
\draw (0.9222029944149257, 1.0535214108139956) -- (1.0557740675891703, 0.9186289404430888);
\draw (0.8378098994600067, 1.0553763756752683) -- (0.9222029944149257, 1.0535214108139956);
\draw (0.5520153085255308, 0.5495884614045772) -- (0.6950107605184371, 0.5157164564996014);
\draw (0.5176560761734916, 0.6886283467274534) -- (0.5520153085255308, 0.5495884614045772);
\draw (0.6641883932565097, 0.8689206838898882) -- (0.5176560761734916, 0.6886283467274534);
\draw (0.7058167568558498, 0.868884105906581) -- (0.6641883932565097, 0.8689206838898882);
\draw (0.8718044537844236, 0.7026819692452801) -- (0.7058167568558498, 0.868884105906581);
\draw (0.8718847063078764, 0.6736083554769277) -- (0.8718044537844236, 0.7026819692452801);
\draw (0.6950107605184371, 0.5157164564996014) -- (0.8718847063078764, 0.6736083554769277);
\draw[dashed] (0, 0) -- (1, 0);
\draw[dashed] (1, 0) -- (1, 1);
\draw[dashed] (1, 1) -- (0, 1);
\draw[dashed] (0, 1) -- (0, 0);
\draw (-0.02, -0.13) node[anchor=north west] { \scriptsize $d_{ 13 }(T^2) \leqslant 0.354547 $ };
\end{tikzpicture} \\
    \end{tabular}
    \caption{Partitions of the torus into $9, 10, 13$ parts.} 
    \label{torus_opt_alg}
\end{figure}
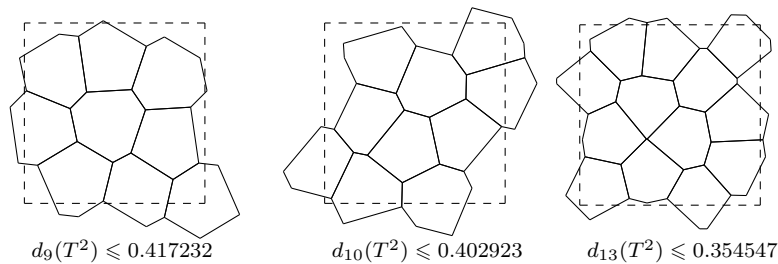

Then the local minimum can be found using a modification of the gradient descent for piecewise smooth functions, such as the basic stochastic gradient descent or the Adam algorithm \cite{adam}. As in previous work \cite{DAM_our_article}, a Voronoi diagram constructed for some random packing of circles is used as an initial approximation (Fig. \ref{alg_example}), and multiple restarts of the optimisation algorithm are performed. Some of the partitions found using this approach are shown in Fig. \ref{torus_opt_alg}.

Note that almost any nonsmooth optimization algorithm can be applied to find the local minimum in \eqref{globopt}. Here we deal with a piecewise quadratic problem without constraints, which from the point of view of optimization methods is simpler than  partitions of plane sets, where constraints were present. Nevertheless, as in the previous work, we used the Adam algorithm because its implementations in machine learning libraries allow us to easily use the capabilities of parallel computations and computations on GPUs. 

 \section{Conclusion}

 As in the problems with packing $m$ circles or spheres in a container or in a compact manifold, it seems that the problem considered cannot be completely solved for all $m$.  However, let us formulate some questions that remain outside the scope of this paper.

 \begin{question}
    Is it true that the estimates of $d_m(T^2)$ at $m = 4,5,6$ given in Theorem \ref{theorem3} are exact? If yes, how to prove it?
 \end{question}

Note that the regions into which the torus is partitioned, in general, do not have to be connected, so it is not clear how to reduce the problem for given $m$ to the finite enumeration.

\begin{question}
    Denote by $M^*(\tau) \subset T^2$ the solution of the following problem.
    \[
        \operatorname{diam} M = \tau,        
    \]
    \[
        \operatorname{Area}(M) \to \max,
    \]
     What is the shape of the set $M^*(\tau)$, if $\frac{1}{2}<\tau<\frac{\sqrt{2}}{2}$?
\end{question}

\begin{figure}[!htb]
    \centering

    \begin{tabular}{cccc}
    \includegraphics[scale=0.19]{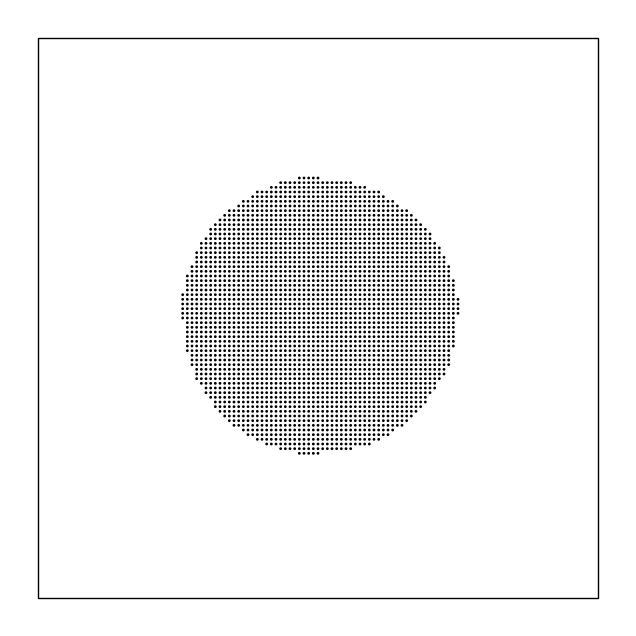} &
    \includegraphics[scale=0.19]{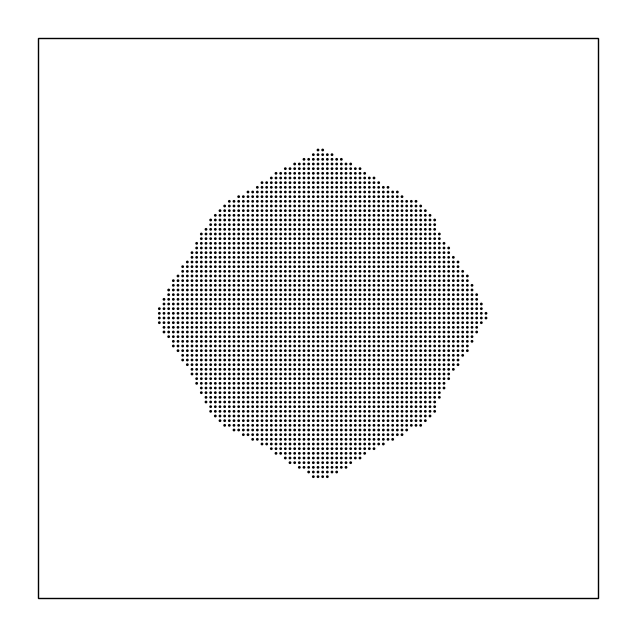} &
    \includegraphics[scale=0.19]{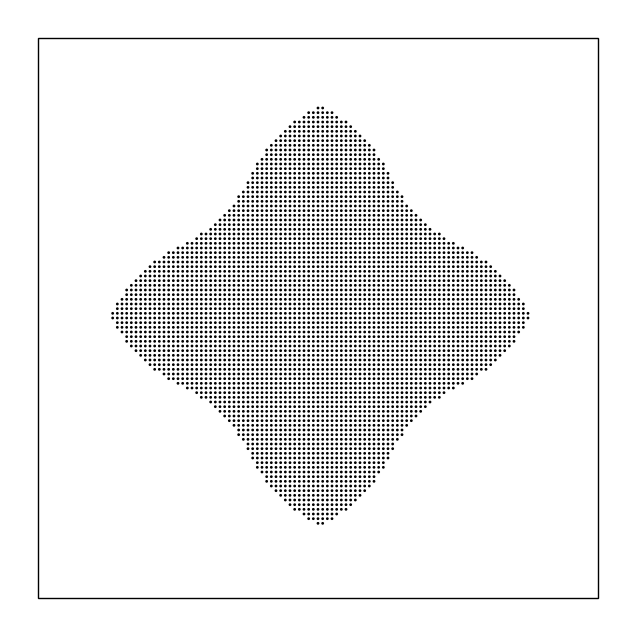} &
    \includegraphics[scale=0.19]{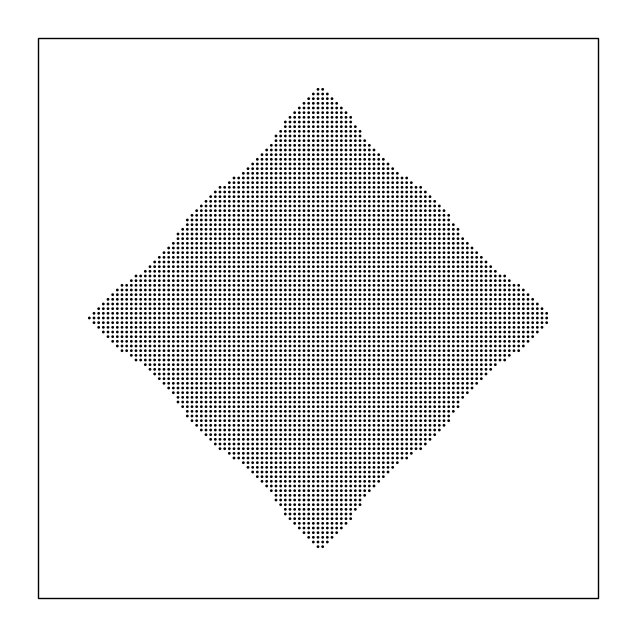} \\
    \end{tabular}
    \caption{Shape of the maximum area, $\tau \in \{ 0.493, 0.527, 0.559, 0.589\}$.} 
    \label{max_area_sets}
\end{figure}

Solving this problem is unlikely to strengthen the estimates of $d_m(T^2)$, but it may be of particular interest. Figure \ref{max_area_sets} represents visualizations of the maximal independent sets found by the KaMIS software \cite{Kamis} on a  $120 \times 120$ grid at different values of $\tau$. 
These figures illustrate the resulting locally optimal configurations that maximize the area within the given constraints in the approximate discrete problem.

\bibliography{torus_partitions}

\end{document}